\documentclass[number,citesort,dvips]{arxbj}
\usepackage{upgreek}

\aid{0}
\volume{17}
\issue{3}
\pubyear{2011}
\firstpage{942}
\lastpage{968}
\doi{10.3150/10-BEJ299}

\makeatletter
\newcommand{\eqref}[1]{(\ref{#1})}
\newcommand{\notag}{\nonumber}

\newcommand{\E}{\mathbb E}
\newcommand{\R}{\mathbb{R}}
\newcommand{\N}{\mathbb{N}}
\newcommand{\Z}{\mathbb{Z}}
\newcommand{\XXX}{\mathbb{X}}
\newcommand{\YYY}{\mathbb{Y}}
\newcommand{\GGG}{\mathbb{G}}
\newcommand{\SSS}{\mathbb{S}}
\newcommand{\afrac}[2]{{#1/#2}}
\newcommand{\bfrac}[2]{{(#1/#2)}}

\renewcommand{\P}{\mathbb{P}}
\newcommand{\Var}{\operatorname{Var}}
\newcommand{\Cov}{\operatorname{Cov}}
\newcommand{\Aa}{\mathcal{A}}
\newcommand{\eps}{\varepsilon}
\newcommand{\tofd}{\stackrel{f.d.d.}{\rightarrow}}
\newcommand{\toweak}{\stackrel{w}{\rightarrow}}
\newcommand{\toskor}{\stackrel{w}{\rightarrow}}
\newtheorem{theorem}{Theorem}[section]
\newtheorem{lemma}{Lemma}[section]
\newtheorem{proposition}{Proposition}[section]
\newremark{remark}{Remark}[section]
\makeatother

\begin{document}
\begin{frontmatter}

\title{Functional limit theorems for sums of independent geometric L\'evy processes}
\runtitle{Sums of independent geometric L\'evy processes}

\begin{aug}
\author{\fnms{Zakhar} \snm{Kabluchko}\corref{}\ead[label=e1]{zakhar.kabluchko@uni-ulm.de}}
\runauthor{Z. Kabluchko}
\address{Institute of Stochastics, Ulm University, Helmholtzstr. 18,
89069 Ulm, Germany}
\end{aug}

\received{\smonth{11} \syear{2009}}
\revised{\smonth{6} \syear{2010}}

%
\begin{abstract}
Let $\xi_i$, $i\in\N$, be independent copies of a L\'evy process $\{
\xi(t), t\geq0\}$.
Motivated by the results obtained previously in the context of the
random energy model, we prove
functional limit theorems for the process
\[
Z_N(t)=\sum_{i=1}^N \mathrm{e}^{\xi_i(s_N+t)}
\]
as $N\to\infty$, where $s_N$ is a non-negative sequence converging to
$+\infty$.
The limiting process depends heavily on the growth rate of the sequence $s_N$.
If $s_N$ grows slowly in the sense that $\liminf_{N\to\infty}\log
N/s_N>\lambda_2$ for some critical value $\lambda_2>0$, then the
limit is an Ornstein--Uhlenbeck process. However, if $\lambda:=\lim
_{N\to\infty}\log N/s_N\in(0,\lambda_2)$, then the limit is a
certain completely asymmetric $\alpha$-stable process $\YYY_{\alpha
;\xi}$.
\end{abstract}

%

\end{frontmatter}

\section{Introduction and statement of main results}\label{sec:intro}
\subsection{Introduction}\label{subsec:intro}
One of the simplest models in the physics of disordered systems is the
random energy model (REM). The partition function of the random energy
model at an inverse temperature $\beta>0$ is a random variable
$S_n(\beta)$ given by
%
\begin{equation}\label{eq:REM}
S_n(\beta)=\sum_{i=1}^{2^n} \mathrm{e}^{\beta\sqrt{n}\zeta_i},
\end{equation}
where $\zeta_i$, $i\in\N$, are i.i.d. standard Gaussian random variables.
Bovier \textit{et al.} \cite{bovieretal02} studied the limit laws of $S_n(\beta)$ as $n\to
\infty$ in dependence on the parameter $\beta$. They showed that for
$\beta<\sqrt{\log2/2}$, the random variable $S_n(\beta)$ obeys a
central limit theorem with a Gaussian limit law, whereas for $\beta
>\sqrt{\log2/2}$, the limit distribution is a completely asymmetric
$\alpha$-stable law. The results of \cite{bovieretal02} have been
extended by Ben Arous \textit{et al.} \cite{benarousetal05} to the case when the random
variables $\zeta_i$ are non-Gaussian; see also \cite
{bogachev06,bogachev07,janssen09}.
Extending~\cite{bovieretal02} in a~different direction, Cranston and
Molchanov~\cite{cranstonmolchanov05} considered sums of the form
%
\begin{equation}\label{eq:REM_CrMol}
R_n(\beta)=\sum_{i=1}^{N(n)} \mathrm{e}^{\beta\sum_{j=1}^n \zeta_{i,j}},
\end{equation}
where $\zeta_{i,j}$, $(i,j)\in\N^2$, is a two-dimensional array of
i.i.d. random variables, $N(n)$ is a~certain exponentially growing
function of $n$, $\beta>0$, and $n\to\infty$. The sum $R_n(\beta)$
reduces to $S_n(\beta)$ if the random variables $\zeta_{i,j}$ are
standard Gaussian and $N(n)=2^n$. Cranston and Molchanov \cite{cranstonmolchanov05} have
shown that the behavior of the sum $R_n(\beta)$ is rather similar to
that of the sum $S_n(\beta)$, with Gaussian and completely asymmetric
$\alpha$-stable limit laws.
Unaware of \cite{cranstonmolchanov05}, the author proved essentially
the same result in \cite{kabluchkoprod}.

The aim of the present paper is to obtain \textit{functional limit
theorems} corresponding to the results of \cite
{bovieretal02,cranstonmolchanov05,kabluchkoprod}. That is, we will
consider sums of exponentials of stochastic processes (L\'evy processes
or random walks) rather than sums of exponentials of random variables.
We prefer to work with L\'evy processes, but it should be stressed that
all our results have straightforward analogues for random walks. Let
$\xi_i$, $i\in\N$, be independent copies of a L\'evy process $\{\xi
(t), t\geq0\}$, and let $\{s_N\}_{N\in\N}$ be a non-negative
sequence. We are interested in the limiting properties, as $N\to\infty
$, of the stochastic process $Z_N$ defined by\looseness=-1
%
\begin{equation}\label{eq:def_w}
Z_N(t)=\sum_{i=1}^{N}\mathrm{e}^{\xi_{i}(s_N+t)}.
\end{equation}

Since the random variable $Z_N(0)$ reduces essentially to
$R_{s_N}(\beta)$, we will recover the results of \cite
{bovieretal02,cranstonmolchanov05,kabluchkoprod} by restricting our
processes to $t=0$. If $s_N=\beta^2 n$, $N=2^n$, and $\xi$ is a
standard Brownian motion, then $Z_N(0)$ has the same distribution as
the partition function of the random energy model $S_n(\beta)$ given
in \eqref{eq:REM}. The results of \cite
{bovieretal02,cranstonmolchanov05,kabluchkoprod} suggest that the
limiting process for $Z_N$ as $N\to\infty$ should be either Gaussian
or completely asymmetric $\alpha$-stable depending on the rate of
growth of the sequence $s_N$. We will show that this is indeed the
case, obtaining in the limit an Ornstein--Uhlenbeck process in the
``slow growth regime'', and a certain completely asymmetric $\alpha
$-stable process $\YYY_{\alpha;\xi}$ in the ``fast growth regime''.
The family of processes $\YYY_{\alpha;\xi}$ has not been studied in
the literature so far, although a similar class of max-stable processes
has been considered in \cite{stoev07}.

To give a motivation for studying the process $Z_N$, consider the
following problem. Suppose that we are given a portfolio consisting of
a large number $N$ of financial assets whose prices are modeled by
independent geometric Brownian motions (or, somewhat more generally, by
independent geometric L\'evy processes). Then, the price of the whole
portfolio after $s_N$ units of time have passed is given by the process
$Z_N$. It will be shown below that if $s_N$, as a function of $N$,
grows slowly (i.e.,  if we are looking at the price in the \textit
{near} future), then the price of the portfolio is approximated by an
Ornstein--Uhlenbeck process, whereas if $s_N$ grows rapidly (i.e., if
we are interested in the \textit{remote} future), then the price is
approximated by the $\alpha$-stable process $\YYY_{\alpha;\xi}$.
For example, if we are summing standard geometric Brownian motions,
then the boundary between the near future and the remote future lies at
$s_N\sim\frac12\log N$.

\subsection{Notation}\label{sec:not_0}
Before we can state our results, we need to recall some facts related
to Cram\'er's large deviations theorem; see, for instance,
\cite{dembobook}, Chapter~2.2. A L\'evy process is a process with stationary,
independent increments and cadlag sample paths. Let $\{\xi(t), t\geq
0\}$ be a L\'evy process such that
%
\begin{equation}\label{eq:def_varphi}
\psi(u):=\log\E \mathrm{e}^{u\xi(1)}  \qquad \mbox{is finite for all } u\in\R.
\end{equation}
We always assume that $\xi(1)$ is not a.s. constant.
The function $\psi$ is infinitely differentiable and strictly convex
with $\psi(0)=0$. It follows that $\psi'\dvtx [0,\infty)\to[\beta
_{0},\beta_{\infty})$ is a monotone increasing bijection, where
%
\begin{equation}\label{eq:def_beta_zero_infty}
\beta_{0}=\psi'(0)=\E\xi(1),   \qquad   \beta_{\infty}=\lim_{u\to
+\infty}\psi'(u).
\end{equation}
Let $I\dvtx [\beta_0,\beta_{\infty})\to[0,+\infty)$ be the
Legendre--Fenchel transform of $\psi$ defined by
%
\begin{equation}\label{eq:I_varphi_prime}
I(\psi'(u))=u\psi'(u)-\psi(u),  \qquad    u\geq0.
\end{equation}
The function $I$ is strictly increasing, strictly convex, infinitely
differentiable with $I(\beta_0)=0$.
As in \cite{cranstonmolchanov05,kabluchkoprod}, it will turn out
that the limiting properties of the process $Z_N$ undergo phase
transitions at the ``critical points'' $\lambda_1,\lambda_2$ given by
%
\begin{equation}\label{eq:def_lambda}
\lambda_1=I(\psi'(1))=\psi'(1)-\psi(1),   \qquad   \lambda_2=I(\psi
'(2))=2\psi'(2)-\psi(2).
\end{equation}
For example, if $\xi$ is a standard Brownian motion, then $\psi
(u)=I(u)=u^2/2$ and the critical points are given by $\lambda_1=1/2$,
$\lambda_2=2$.

\subsection{Statement of main results}\label{subsec:res}
Our first result deals with the case $s_N=0$ (but covers automatically
also the case $s_N=\mathit{const}$). It is a consequence of the central limit
theorem in the Skorokhod space, and is stated merely for completeness.
%
\begin{theorem}\label{theo:zero}
If $s_N=0$ and condition \eqref{eq:def_varphi} holds, then for every
$T>0$, we have the following weak convergence of stochastic processes
on the Skorokhod space $D[0,T]$:
%
\begin{equation}\label{eq:theo_0}
\frac{Z_N(\cdot)-\E Z_N(\cdot)}{\sqrt{N}}\toskor
\GGG(\cdot),  \qquad    N\to\infty,
\end{equation}
where $\{\GGG(t), t\geq0\}$ is a zero-mean Gaussian process with
covariance function
%
\begin{equation}\label{eq:cov_G}
\Cov(\GGG(t_1),\GGG(t_2))=\mathrm{e}^{\psi(2)t_1+\psi(1)(t_2-t_1)}-\mathrm{e}^{\psi
(1)(t_1+t_2)}, \qquad
0\leq t_1\leq t_2.
\end{equation}
\end{theorem}

Our next theorem deals with the case in which $s_N$ grows slowly as a
function of $N$. We will assume that the following slow growth
condition is satisfied:
%
\begin{equation}\label{eq:asymps_N1}
\lim_{N\to\infty}s_N=\infty, \qquad    \liminf_{N\to\infty} \frac
{\log N}{s_N}>\lambda_2.
\end{equation}
%
%
\begin{theorem}\label{theo:normal}
If conditions \eqref{eq:def_varphi} and \eqref{eq:asymps_N1} hold,
then for every $T>0$, we have the following weak convergence of
stochastic processes on the Skorokhod space $D[-T,T]$:
%
\begin{equation}\label{eq:theo_normal}
\frac{Z_{N}(\cdot)-\E Z_{N}(\cdot)}{\sqrt{\Var Z_N(\cdot)}}
\toskor\XXX(\cdot),   \qquad    N\to\infty,
\end{equation}
where
$\{\XXX(t), t\in\R\}$ is a zero-mean Gaussian process with
covariance function
%
\begin{equation}\label{eq:cov_X}
\Cov(\XXX(t_1),\XXX(t_2))=\mathrm{e}^{(\psi(1)-\afrac{\psi
(2)}{2})|t_2-t_1|},  \qquad   t_1,t_2\in\R.
\end{equation}
\end{theorem}

Note that $\XXX$ is an Ornstein--Uhlenbeck process and that the
process on the left-hand side of \eqref{eq:theo_normal} is
well defined on $[-T,T]$ if $N$ is sufficiently large.
In the next theorem, which deals with the ``critical case'', we still
obtain an Ornstein--Uhlenbeck process in the limit, but an additional
factor appears. We will assume that the following critical growth
condition holds: For some $\vartheta\in\R$,
%
\begin{equation}\label{eq:asymps_Ncrit}
\log N=\lambda_2
s_N+2\vartheta\sqrt{\psi''(2)s_N}+\mathrm{o}\bigl(\sqrt{s_N}\bigr),\qquad
   N\to\infty.
\end{equation}
\begin{theorem}\label{theo:crit}
If conditions \eqref{eq:def_varphi} and \eqref{eq:asymps_Ncrit} are
satisfied, then we have the following convergence of stochastic processes:
%
\begin{equation}
\frac{Z_N(\cdot)-\E Z_N(\cdot)}{\sqrt{\Var Z_N(\cdot)}}
\tofd\sqrt{\Phi(\vartheta)} \XXX(\cdot),  \qquad    N\to\infty,
\end{equation}
where $\Phi$ is the standard normal distribution function, $\XXX$ is
as in Theorem \ref{theo:normal}, and $\tofd$ denotes the weak
convergence of finite-dimensional distributions.
\end{theorem}

Let us stress that even when restricted to $t=0$, the above theorem
gives a more ``smooth'' picture of the critical regime than the
corresponding results of \cite{bovieretal02,cranstonmolchanov05,kabluchkoprod} where only the case $\vartheta=0$
has been considered.

The next theorem shows that in the fast growth case, a non-Gaussian
process $\YYY_{\alpha;\xi}$ appears in the limit.
We need the following fast growth condition:
%
\begin{equation}\label{eq:asymps_N}
\lambda:=\lim_{N\to\infty}\frac{\log N}{s_N}\in(0,\lambda_2).
\end{equation}
Recall also that a random variable is called lattice if its values are
of the form $an+b$, $n\in\Z$, for some $a,b\in\R$, and non-lattice
if no such $a$ and $b$ exist.
\begin{theorem}\label{theo:stable}
Suppose that \eqref{eq:def_varphi} and \eqref{eq:asymps_N} hold, and
assume that the distribution of $\xi(1)$ is non-lattice. Define
$\alpha\in(0,2)$ as the unique solution of the equation $I(\psi
'(\alpha))=\lambda$ and~let
%
\begin{equation}\label{eq:def_An}
A_N(t)
=
\cases{\displaystyle
0, &\quad  if  $\lambda\in(0,\lambda_1)$,\cr\displaystyle
\mathrm{e}^{\psi(1)t} N \E\bigl[\mathrm{e}^{\xi(s_N)}1_{\xi(s_N)\leq\log
B_N(0)}\bigr]+l(t)B_N(t), &\quad  if   $\lambda=\lambda_1$,\cr\displaystyle
\mathrm{e}^{\psi(1)t} \E Z_N(0), &\quad   if  $\lambda\in(\lambda
_1,\lambda_2)$,
}
\end{equation}
where $l(t)=(\psi'(0)-\psi'(1))t1_{t<0}$, and
%
\begin{equation}\label{eq:def_Bn_large}
B_N(t)
=\mathrm{e}^{\bfrac{\psi(\alpha)}{\alpha}t}\exp \biggl\{s_NI^{-1} \biggl(\frac
{\log N-\log( \alpha\sqrt{2\uppi\psi''(\alpha) s_N}) }{s_N}
\biggr) \biggr\}.
\end{equation}
Then, for every $T>0$, we have the following convergence of stochastic
processes on the Skorokhod space $D[-T,T]$:
%
\begin{equation}\label{eq:theo_stable}
\frac{Z_N(\cdot)-A_N(\cdot)}{B_N(\cdot)}\toskor\YYY_{\alpha;\xi
}(\cdot),   \qquad    N\to\infty.
\end{equation}
Here, $\YYY_{\alpha;\xi}$ is a completely asymmetric $\alpha
$-stable process that will be defined below.
\end{theorem}

\begin{remark}
Our results have straightforward discrete-time analogues with geometric
L\'evy processes replaced by exponentials of independent random walks.
If $\xi$ is the standard Brownian motion, then in all our results the
weak convergence in the Skorokhod space can be replaced by the weak
convergence in the space of continuous functions. The non-lattice
assumption in Theorem \ref{theo:stable} cannot be dropped; see \cite
{kabluchkoprodlatt}.
\end{remark}

\subsection{Definition of the process $\YYY_{\alpha;\xi}$}\label
{subsec:def_proc_Y_axi}

We now define the $\alpha$-stable process $\YYY_{\alpha;\xi}$ which
appeared in Theorem \ref{theo:stable}. Our main reference on $\alpha
$-stable distributions and processes is \cite{samorodnitskytaqqubook}.
First of all, fix some $\alpha\in(0,2)$, and let $\xi_i$, $i\in\N
$, be independent copies of a L\'evy process $\{\xi(t), t\geq0\}$
satisfying condition \eqref{eq:def_varphi}.
Independently, let $\{\Gamma_i, i\in\N\}$ be the arrivals of a unit
intensity Poisson process on the positive half-line. In other words,
$\Gamma_k=\sum_{i=1}^k \eps_i$, where $\eps_i$, $i\in\N$, are
i.i.d. exponential random variables with mean $1$. Define $U_i=\Gamma
_i^{-1/\alpha}$, $i\in\N$, and note that $\{U_i,i\in\N\}$ are the
points of a Poisson process on $(0,\infty)$ with intensity $\alpha
u^{-(\alpha+1)}\,\mathrm{d}u$, arranged in the descending order.
The restriction of the process $\YYY_{\alpha;\xi}$ to the positive
half-line is then defined as follows: For $t\geq0$, we set
%
\begin{equation}\label{eq:def_Y}
\YYY_{\alpha;\xi}(t)=
\cases{\displaystyle
\sum_{i\in\N} U_i \mathrm{e}^{\xi_i(t)-\bfrac{\psi(\alpha)}{\alpha}t},
&\quad$0<\alpha<1$,\cr\displaystyle
\lim_{\tau\downarrow0} \biggl(\mathop{\mathop{\sum}_{i\in\N}}_{ U_i>\tau}
U_i \mathrm{e}^{\xi_i(t)-\psi(1)t}-\log\frac{1}{\tau} \biggr),
&\quad$\alpha=1$,\cr\displaystyle
\lim_{\tau\downarrow0}
 \biggl(\mathop{\mathop{\sum}_{i\in\N}}_{ U_i>\tau} U_i \mathrm{e}^{\xi_i(t)-\bfrac
{\psi(\alpha)}{\alpha}t}
-\frac{\alpha \tau^{1-\alpha}}{\alpha-1} \mathrm{e}^{(\psi(1)-\afrac{\psi
(\alpha)}{\alpha})t} \biggr),
&\quad$1<\alpha<2$.
}
\end{equation}
For the definition of the process $\YYY_{\alpha;\xi}$ on the
negative half-line we refer to \cite{kabluchkoprep}.
The Poisson representation of $\alpha$-stable random vectors --
see \cite{samorodnitskytaqqubook}, Theorem 3.12.2 -- implies that
for every $t\geq0$, the expression defining $\YYY_{\alpha;\xi}(t)$
converges with probability $1$. Further, the finite-dimensional
distributions of the process $\YYY_{\alpha;\xi}$ are $\alpha
$-stable with skewness parameter $\beta=1$.
If $\alpha\in(0,1)$, then the process $\YYY_{\alpha;\xi}$ takes
only positive values; otherwise, it takes any real values.
For the proof of the next proposition we refer to \cite{kabluchkoprep}.
\begin{proposition}\label{prop:uniform}
The expression on the right-hand side of \eqref{eq:def_Y} defining
$\YYY_{\alpha;\xi}$ converges uniformly on compact sets with
probability $1$.
\end{proposition}

As a consequence, the process $\YYY_{\alpha;\xi}$ has cadlag sample
paths. Moreover, if $\xi$ is a Brownian motion, then the sample paths
of $\YYY_{\alpha;\xi}$ are even continuous. The process $\YYY
_{\alpha;\xi}$ is stationary for $\alpha\neq1$; see the preprint
version of this paper \cite{kabluchkoprep} for this and other
properties of $\YYY_{\alpha;\xi}$. The rest of the paper is devoted
to proofs.


\section{Large deviations and truncated exponential moments}
The next proposition on the asymptotic behavior of truncated
exponential moments will play a crucial role in the sequel. Parts of it
are scattered over \cite{cranstonmolchanov05,kabluchkoprod}, but we
will give a~simple unified proof below.
\begin{proposition}\label{prop:trunc_clt}
Let $\{\xi(t), t\geq0\}$ be a L\'evy process satisfying \eqref
{eq:def_varphi} and suppose that the distribution of $\xi(1)$ is
non-lattice. Let $\kappa\geq0$, and let $b_N\to\infty$ and $x_N\to
\infty$ be two sequences. Let $I$ be the large deviation function of
$\xi(1)$, as defined in \eqref{eq:I_varphi_prime}.
\begin{enumerate}[(3)]
\item[(1)]\hypertarget{p:1_trunc}
If for some $\vartheta\in\R$, $b_N=\psi'(\kappa)x_N+\vartheta
\sqrt{\psi''(\kappa)x_N}+\mathrm{o}(\sqrt{x_N})$ as $N\to\infty$, then
%
\begin{equation}\label{eq:trunc_p1}
\lim_{N\to\infty} \mathrm{e}^{-\psi(\kappa) x_N} \E\bigl[\mathrm{e}^{\kappa\xi(x_N)}
1_{\xi(x_N)\leq b_N}\bigr]=\Phi(\vartheta),
\end{equation}
where $\Phi$ is the standard Gaussian distribution function.
\item[(2)]\hypertarget{p:2_trunc}
If $\liminf_{N\to\infty} b_N/x_N>\psi'(\kappa)$, then
%
\begin{equation}\label{eq:trunc_p2}
\lim_{N\to\infty} \mathrm{e}^{-\psi(\kappa) x_N} \E\bigl[\mathrm{e}^{\kappa\xi(x_N)}
1_{\xi(x_N)> b_N}\bigr]=0.
\end{equation}
If, moreover, $\lim_{N\to\infty}b_N/x_N=\psi'(\alpha)$ for some
$\alpha>\kappa$, then
%
\begin{equation}\label{eq:trunc_p2a}
\E\bigl[\mathrm{e}^{\kappa\xi(x_N)} 1_{\xi(x_N)> b_N}\bigr]\sim\frac{\mathrm{e}^{\kappa
b_N}}{(\alpha-\kappa) \sqrt{2\uppi\psi''(\alpha)x_N}}
\mathrm{e}^{-I(b_N/x_N)x_N},    \qquad  N\to\infty.
\end{equation}
\item[(3)]\hypertarget{p:3_trunc}
If $\limsup_{N\to\infty} b_N/x_N<\psi'(\kappa)$, then
%
\begin{equation}\label{eq:trunc_p3}
\lim_{N\to\infty} \mathrm{e}^{-\psi(\kappa) x_N}\E\bigl[\mathrm{e}^{\kappa\xi(x_N)}
1_{\xi(x_N)\leq b_N}\bigr]=0.
\end{equation}
If, moreover, $\lim_{N\to\infty}b_N/x_N=\psi'(\alpha)$ for some
$\alpha\in(0,\kappa)$, then
%
\begin{equation}\label{eq:trunc_p3a}
\E\bigl[\mathrm{e}^{\kappa\xi(x_N)} 1_{\xi(x_N)\leq b_N}\bigr]\sim\frac{\mathrm{e}^{\kappa
b_N}}{(\kappa-\alpha) \sqrt{2\uppi\psi''(\alpha)x_N}}
\mathrm{e}^{-I(b_N/x_N)x_N},  \qquad    N\to\infty.
\end{equation}
\end{enumerate}
\end{proposition}

The following precise form of Cram\'er's large deviations theorem was
stated and proved in \cite{bahadurrao60,petrov65} for sums of i.i.d.
random variables, but applies equally well to L\'evy processes.
\begin{theorem}\label{theo:ld}
Let $\{\xi(t), t\geq0\}$ be a L\'evy process satisfying \eqref
{eq:def_varphi} and suppose that the distribution of $\xi(1)$ is non-lattice.
Let $\beta=\psi'(\alpha)$, where $\alpha>0$. Then,
%
\begin{equation}
\P[\xi(T)\geq\beta T]\sim\frac{1}{\alpha\sqrt{2\uppi\psi
''(\alpha)T}} \mathrm{e}^{- I(\beta)T}, \qquad
    T\to\infty.
\end{equation}
The statement holds uniformly in $\beta\in K$ for any compact set
$K\subset(\beta_0,\beta_{\infty})$.
\end{theorem}
\begin{pf*}{Proof of Proposition \ref{prop:trunc_clt}}
We will use an exponential change of measure argument. Denote by $F_t$
the distribution function of $\xi(t)$. There exists a L\'evy process
$\{\tilde\xi(t), t\geq0\}$ (an exponential twist of $\xi$) such
that $\tilde F_t$, the distribution function of $\tilde\xi(t)$, is
given by
%
\begin{equation}\label{eq:tilde_F_F}
\frac{\tilde F_t(\mathrm{d}x)}{F_t(\mathrm{d}x)}=\mathrm{e}^{\kappa x-\psi(\kappa)t},
\qquad
x\in\R.
\end{equation}
Recall from \eqref{eq:def_varphi} that $\psi(u)=\log\E \mathrm{e}^{u\xi(1)}$
and let $\tilde\psi(u)=\log\E \mathrm{e}^{u\tilde\xi(1)}$. By \eqref
{eq:tilde_F_F}, we have
%
\begin{equation}\label{eq:psi_tilde_psi}
\tilde\psi(u)=\log\int_{\R} \mathrm{e}^{ux}\,\mathrm{d}\tilde F_1(x)=\log\int_{\R}
\mathrm{e}^{ux}\mathrm{e}^{\kappa x-\psi(\kappa)}\,\mathrm{d}F_1(x)=\psi(u+\kappa)-\psi(\kappa).
\end{equation}
%
Hence,
%
\begin{equation}\label{eq:E_tilde X}
\E\tilde\xi(T)=\tilde\psi'(0)T=\psi'(\kappa)T, \qquad
\Var\tilde\xi(T) = \tilde\psi''(0)T= \psi''(\kappa)T.
\end{equation}
%
The study of the truncated exponential moment
%
\begin{equation}
M_N:=\mathrm{e}^{-\psi(\kappa) x_N}\E\bigl[\mathrm{e}^{\kappa\xi(x_N)} 1_{\xi(x_N)\leq b_N}\bigr]
\end{equation}
can be reduced to the study of the probability $\P[\tilde\xi
(x_N)\leq b_N]$ as follows:
%
\begin{equation}\label{eq:exp_change}
M_N
=\int_{-\infty}^{b_N}\mathrm{e}^{\kappa x-\psi(\kappa)x_N}\,\mathrm{d}F_{x_N}(x)
=\int_{-\infty}^{b_N}\,\mathrm{d}\tilde F_{x_N}(x) 
=\P[\tilde\xi(x_N)\leq b_N]. 
\end{equation}
%
Having the central limit theorem in mind, we write
%
\begin{equation}\label{eq:S_n_Phi_n}
\P[\tilde\xi(x_N)\leq b_N]=\P \biggl[\frac{\tilde\xi(x_N)-\psi
'(\kappa) x_N}{\sqrt{\psi''(\kappa) x_N}}\leq
r_N \biggr],  \qquad \mbox{where } r_N=\frac{b_N-\psi'(\kappa)
x_N}{\sqrt{\psi''(\kappa) x_N}}.
\end{equation}

Let us prove part \hyperlink{p:1_trunc}{1} of the proposition. By the
assumption of part \hyperlink{p:1_trunc}{1}, we have $\lim_{N\to\infty
}r_N=\vartheta$. Then, it follows from \eqref{eq:exp_change} and the
central limit theorem that
\[
\lim_{N\to\infty}M_N=\lim_{N\to\infty} \P[\tilde\xi(x_N)\leq
b_N]=\Phi(\vartheta),
\]
which proves \eqref{eq:trunc_p1}.

Let us prove part \hyperlink{p:2_trunc}{2}  of the proposition. If $\liminf
_{N\to\infty} b_N/x_N>\psi'(\kappa)$, then $\lim_{N\to\infty
}r_N=+\infty$, and the central limit theorem implies that
\[
\lim_{N\to\infty}M_N=\lim_{N\to\infty} \P[\tilde\xi(x_N)\leq b_N]=1,
\]
which proves \eqref{eq:trunc_p2}.
To prove \eqref{eq:trunc_p2a}, we will apply Theorem \ref{theo:ld} to
the process $\tilde\xi$. The large deviation function of the process
$\tilde\xi$ is defined by $\tilde I(\tilde\psi'(u))=u\tilde\psi
'(u)-\tilde\psi(u)$. Hence, setting $\beta=\tilde\psi'(u)$ and
taking into account \eqref{eq:psi_tilde_psi}, we obtain
\[
\tilde I(\beta)=\tilde I(\tilde\psi'(u))=u\tilde\psi'(u)-\tilde
\psi(u)=u\psi'(u+\kappa)-\psi(u+\kappa)+\psi(\kappa).
\]
Note that $\beta=\psi'(u+\kappa)$ by \eqref{eq:psi_tilde_psi}. It
follows that we have the following formula for the function~$\tilde I$:
%
\begin{equation}\label{eq:def_tilde_I}
\tilde I(\beta)= I(\beta)+\psi(\kappa)-\kappa\beta.
\end{equation}
If $\lim_{N\to\infty}b_N/x_N=\psi'(\alpha)=\tilde\psi'(\alpha
-\kappa)$, then we apply Theorem \ref{theo:ld} to obtain that
\[
\P[\tilde\xi(x_N)>b_N]\sim\frac{1}{(\alpha-\kappa) \sqrt{2\uppi
\psi''(\alpha)x_N}} \mathrm{e}^{- \tilde I(b_N/x_N)x_N},   \qquad   N\to\infty.
\]
A straightforward calculation using \eqref{eq:def_tilde_I} leads
to \eqref{eq:trunc_p2a}.
The proof of part \hyperlink{p:3_trunc}{3} of the proposition is analogous to
the proof of part \hyperlink{p:2_trunc}{2}.
\end{pf*}

We will need the following lemmas; see \cite{kabluchkoprod}, Lemma 3,
and \cite{kabluchkoprep}, Lemma 8.1, for their proofs.
\begin{lemma}\label{lem:I_prime_psi_prime}
For every $u>0$, $I'(\psi'(u))=u$.
\end{lemma}
\begin{lemma}\label{lem:small_incr}
Let $\xi$ be a L\'evy process satisfying \eqref{eq:def_varphi}. Let
$p\in[1,2]$ and fix some $T>0$. Then, there is $C>0$ such that for all
$t\in[0,T]$,
%
\begin{equation}\label{eq:small_incr_st}
\E\bigl|\mathrm{e}^{\xi(t)}-1\bigr|^{p}\leq C t^{p/2},  \qquad   \E\bigl|\mathrm{e}^{2\xi(t)}-\mathrm{e}^{\xi
(t)}\bigr|^{p}\leq C t^{p/2}.
\end{equation}
\end{lemma}

\section{\texorpdfstring{Proof of Theorem \protect\ref{theo:zero}}{Proof of Theorem 1.1}}

The proof is a standard application of the central limit theorem in the
Skorokhod space. First let us compute the covariance function of the
process $\mathrm{e}^{\xi}$. We have, for $0\leq t_1\leq t_2$,
\[
\E\bigl[\mathrm{e}^{\xi(t_1)}\mathrm{e}^{\xi(t_2)}\bigr]=\E \mathrm{e}^{2\xi(t_1)}\cdot\E \mathrm{e}^{\xi
(t_2)-\xi(t_1)}
=\mathrm{e}^{\psi(2)t_1+\psi(1)(t_2-t_1)}.
\]
Since $\E \mathrm{e}^{\xi(t)}=\mathrm{e}^{\psi(1)t}$, we have
\[
\Cov\bigl(\mathrm{e}^{\xi(t_1)},\mathrm{e}^{\xi(t_2)}\bigr)=\mathrm{e}^{\psi(2)t_1+\psi
(1)(t_2-t_1)}-\mathrm{e}^{\psi(1)(t_1+t_2)}.
\]
An application of the multidimensional central limit theorem proves
that \eqref{eq:theo_0} holds in the sense of the weak convergence of
finite-dimensional distributions. To prove the weak convergence in the
space $D[0,T]$, we will verify the conditions of \cite{hahn78}, Theorem 2.
For every $0\leq t_1\leq t_2\leq T$, we have
\[
\E\bigl(\mathrm{e}^{\xi(t_2)}-\mathrm{e}^{\xi(t_1)}\bigr)^2=
\E \mathrm{e}^{2\xi(t_1)}\cdot\E\bigl(\mathrm{e}^{\xi(t_2)-\xi(t_1)}-1\bigr)^2<C(t_2-t_1),
\]
where the last inequality follows from Lemma \ref{lem:small_incr}.
This verifies the first condition of \cite{hahn78}, Theorem 2. The
second condition can be proved in a similar way: for every $0\leq
t_1\leq t_2\leq t_3\leq T$, we have
\begin{eqnarray*}
&&\E\bigl[\bigl(\mathrm{e}^{\xi(t_2)}-\mathrm{e}^{\xi(t_1)}\bigr)^2\bigl(\mathrm{e}^{\xi(t_3)}-\mathrm{e}^{\xi
(t_2)}\bigr)^2\bigr]\\
&& \quad =\E \mathrm{e}^{2\xi(t_1)} \cdot\E\bigl(\mathrm{e}^{2(\xi(t_2)-\xi(t_1))}-\mathrm{e}^{\xi
(t_2)-\xi(t_1)}\bigr)^2\cdot\E\bigl(\mathrm{e}^{\xi(t_3)-\xi(t_2)}-1\bigr)^2\\
&& \quad =\E \mathrm{e}^{2\xi(t_1)}\cdot\E\bigl(\mathrm{e}^{2\xi(t_2-t_1)}-\mathrm{e}^{\xi(t_2-t_1)}\bigr)^2
\cdot\E\bigl(\mathrm{e}^{\xi(t_3-t_2)}-1\bigr)^2\\
&& \quad \leq C(t_3-t_1)^2,
\end{eqnarray*}
where the last inequality follows from Lemma \ref{lem:small_incr}.
This completes the proof.

\section{\texorpdfstring{Proof of Theorem \protect\ref{theo:normal}}{Proof of Theorem 1.2}}
\label{sec:proof_normal}
\subsection{Weak convergence of finite-dimensional distributions}
The first step in establishing Theorem \ref{theo:normal} is to prove
the weak convergence of finite-dimensional distributions in \eqref
{eq:theo_normal}.
It will be convenient to define a positive-valued stochastic process
$W_{N}$ by
%
\begin{equation}\label{eq:def_bn}
W_{N}(t)=N^{-1/2}\mathrm{e}^{\xi(s_N+t)-\bfrac{\psi(2)}{2}(s_N+t)}.
\end{equation}
Let $t_1\leq\cdots \leq t_d$ be fixed, and define a $d$-dimensional
random vector $\mathbf{W}_{N}=(W_{N}(t_1),\ldots ,\break  W_{N}(t_d))$.
If $\mathbf{W}_{1,N},\ldots ,\mathbf{W}_{N,N}$ are independent copies
of $\mathbf{W}_N$, then our aim is to prove that
%
\begin{equation}\label{eq:theo_normal_eq_fd}
\sum_{i=1}^{N} (\mathbf{W}_{i,N}- \E\mathbf{W}_{i,N} )
\toweak
(\XXX(t_k))_{k=1}^d, \qquad      N\to\infty.
\end{equation}
To see that this implies the weak convergence of finite-dimensional
distributions in Theorem~\ref{theo:normal}, it suffices to show that
$\Var Z_N(t) \sim N\mathrm{e}^{\psi(2)(s_N+t)}$ as $N\to\infty$. This can be
done as follows:
%
\begin{eqnarray}\label{eq:var_asympt}
\Var Z_N(t)
&=&N\bigl(\E \mathrm{e}^{2\xi(s_N+t)}-\bigl(\E \mathrm{e}^{\xi(s_N+t)}\bigr)^2\bigr)\notag\\
&=&N\bigl(\mathrm{e}^{\psi(2)(s_N+t)}-\mathrm{e}^{2\psi(1)(s_N+t)}\bigr)\\
&\sim& N \mathrm{e}^{\psi(2)(s_N+t)},  \qquad    N\to\infty,\notag
\end{eqnarray}
where we have used that $\lim_{N\to\infty}s_N=\infty$ by \eqref
{eq:asymps_N1} and that $\psi(2)>2\psi(1)$ by the strict convexity of
$\psi$.


We start proving \eqref{eq:theo_normal_eq_fd}. First of all, let us
compute the covariance matrix of the random vector $\mathbf{W}_N$.
Using \eqref{eq:def_bn} and \eqref{eq:def_varphi}, as well as the
fact that $\xi$ is a L\'evy process, we obtain that for every $1\leq
k\leq l\leq d$,
%
\begin{eqnarray}\label{eq:covar_WN}
\E[W_{N}(t_k)W_{N}(t_l)]
&=& N^{-1}\mathrm{e}^{-\psi(2)s_N}\mathrm{e}^{-\bfrac{\psi(2)}{2}(t_k+t_l)} \E \mathrm{e}^{\xi
(s_N+t_k)+\xi(s_N+t_l)}\notag \\
&=& N^{-1}\mathrm{e}^{-\psi(2)s_N}\mathrm{e}^{-\bfrac{\psi(2)}{2}(t_k+t_l)}
\E \mathrm{e}^{2\xi(s_N+t_k)} \cdot\E
\mathrm{e}^{\xi(s_N+t_l)-\xi(s_N+t_k)}\nonumber
\\[-8pt]
\\[-8pt]
&=& N^{-1}\mathrm{e}^{-\psi(2)s_N}\mathrm{e}^{-\bfrac{\psi(2)}{2}(t_k+t_l)}
\mathrm{e}^{\psi(2)(s_N+t_k)} \mathrm{e}^{\psi(1)(t_l-t_k)}\notag\\
&=& N^{-1} \mathrm{e}^{(\psi(1)-\afrac{\psi(2)}{2})(t_l-t_k)}.\notag
\end{eqnarray}
Since $\psi(2)>2\psi(1)$ by the strict convexity of $\psi$, and
$\lim_{N\to\infty}s_N=\infty$ by \eqref{eq:asymps_N1}, we have for
every $k=1,\ldots ,d$,
%
\begin{equation}\label{eq:sqrtN_EW}
\sqrt{N} \E W_N(t_k)=\mathrm{e}^{\psi(1)(s_N+t_k)}\mathrm{e}^{-\bfrac{\psi
(2)}{2}(s_N+t_k)}\to0,   \qquad   N\to\infty.
\end{equation}
It follows from \eqref{eq:covar_WN} and \eqref{eq:sqrtN_EW} that
%
\begin{equation}\label{eq:covar_WN1}
\lim_{N\to\infty}N\Cov(W_{N}(t_k),W_{N}(t_l))=\mathrm{e}^{(\psi(1)-\afrac
{\psi(2)}{2})(t_l-t_k)}=\Cov(\XXX(t_k),\XXX(t_l)).
\end{equation}

In order to establish \eqref{eq:theo_normal_eq_fd}, we will verify the
Lindeberg condition, that is, we will show that for every $\eps>0$,
%
\begin{equation}\label{eq:lindeberg}
\lim_{N\to\infty} N \E\bigl[\|\mathbf{W}_N-\E\mathbf{W}_N\|^21_{\|
\mathbf{W}_N-\E\mathbf{W}_N\|>\eps}\bigr]=0,
\end{equation}
where $\|\cdot\|$ is the Euclidean norm on $\R^d$.
The multivariate form of the Lindeberg condition we are using can be
found, for example, in \cite{araujoginebook}, Example 4 on page~41.
Since $\lim_{N\to\infty}\sqrt{N} \E\mathbf{W}_N=0$ by \eqref
{eq:sqrtN_EW}, we have $\|\E\mathbf{W}_N\|<\eps/2$ for $N$ large
enough. Thus, for $N$ large enough,
%
\begin{equation}\label{eq:lindeberg_wspom}
\E\bigl[\|\mathbf{W}_N-\E\mathbf{W}_N\|^21_{\|\mathbf{W}_N-\E\mathbf
{W}_N\|>\eps}\bigr]
\leq\E\bigl[\|\mathbf{W}_N-\E\mathbf{W}_N\|^2 1_{\|\mathbf{W}_N\|>\eps/2}\bigr].
\end{equation}
Applying the inequality $\|w_1+w_2\|^2\leq2\|w_1\|^2+2\|w_2\|^2$ to
the right-hand side of \eqref{eq:lindeberg_wspom}, we get
\[
N\E\bigl[\|\mathbf{W}_N-\E\mathbf{W}_N\|^21_{\|\mathbf{W}_N-\E\mathbf
{W}_N\|>\eps}\bigr]
\leq2N\E\bigl[\|\mathbf{W}_N\|^2 1_{\|\mathbf{W}_N\|>\eps/2}\bigr]+2N\|\E
\mathbf{W}_N\|^2.
\]
Note that the second term on the right-hand side converges to $0$
by \eqref{eq:sqrtN_EW}. Hence, in order to prove \eqref
{eq:lindeberg}, it suffices to show that for every $\eps>0$,
%
\begin{equation}\label{eq:lindeberg1}
\lim_{N\to\infty} N\E\bigl[\|\mathbf{W}_N\|^2 1_{\|\mathbf{W}_N\|>\eps}\bigr]=0.
\end{equation}
Let $\Aa_{N,k}$, $k=1,\ldots ,d$, be the random event $\{W_N(t_k)\geq
W_N(t_l), l=1,\ldots ,d\}$. On $\Aa_{N,k}$, we have $\|\mathbf{W}_N\|
^2\leq d W_N^2(t_k)$. Hence,
\begin{eqnarray*}
\E\bigl[\|\mathbf{W}_N\|^2 1_{\|\mathbf{W}_N\|>\eps}\bigr]
&\leq&\sum_{k=1}^d \E\bigl[\|\mathbf{W}_N\|^2 1_{\|\mathbf{W}_N\|>\eps
}1_{\Aa_{N,k}}\bigr]\\
&\leq& d\sum_{k=1}^d \E\bigl[W_N^2(t_k) 1_{W_N(t_k)>\eps/\sqrt d}\bigr].
\end{eqnarray*}
Thus, in order to prove \eqref{eq:lindeberg}, it suffices to show that
for every $t\in\R$ and every $\eps>0$,
%
\begin{equation}\label{eq:clt_wspom1}
\lim_{N\to\infty} N\E\bigl[W_N^2(t) 1_{W_N(t)>\eps}\bigr]=0.
\end{equation}
Recalling \eqref{eq:def_bn} and setting $x_N=s_N+t$ and $b_N=\frac12
(\log N+\psi(2)x_N)+\log\eps$, we may write
%
\begin{equation}\label{eq:wspom555}
N\E\bigl[W_N^2(t) 1_{W_N(t)>\eps}\bigr]=\mathrm{e}^{-\psi(2)x_N}\E\bigl[\mathrm{e}^{2\xi
(x_N)}1_{\xi(x_N)> b_N}\bigr].
\end{equation}
Note that by the slow growth condition \eqref{eq:asymps_N1},
\[
\liminf_{N\to\infty} \frac{b_N}{x_N}>\frac12\bigl(\lambda_2+\psi
(2)\bigr)=\psi'(2).
\]
Applying part \hyperlink{p:2_trunc}{2} of Proposition \ref{prop:trunc_clt}
with $\kappa=2$ to the right-hand side of \eqref{eq:wspom555} we
obtain~\eqref{eq:clt_wspom1}. This verifies the Lindeberg
condition \eqref{eq:lindeberg} and, together with \eqref
{eq:covar_WN1}, completes the proof of the weak convergence of
finite-dimensional distributions in Theorem \ref{theo:normal}.

\subsection{Tightness}
In the rest of the section we complete the proof of Theorem \ref
{theo:normal} by showing that the sequence
%
\begin{equation}\label{eq:tight_seq}
 \biggl\{\frac{Z_N(t)-\E Z_N(t)}{\sqrt{\Var Z_N(t)}}, t\in
[-T,T] \biggr\}_{N\in\N}
\end{equation}
is a tight sequence of stochastic processes in the Skorokhod space
$D[-T,T]$, where $T>0$ is fixed. Since the sequence \eqref
{eq:tight_seq} does not change if we replace the L\'evy process $\xi$
by the L\'evy process $\tilde\xi(t):=\xi(t)-\psi(1)t$, we may and
will assume that
%
\begin{equation}\label{eq:tight_simpl_ass}
\E \mathrm{e}^{\xi(t)}=1,     \qquad  t\geq0.
\end{equation}
Further, since by \eqref{eq:var_asympt}, $\Var Z_N(t)\sim N \mathrm{e}^{\psi
(2)(s_N+t)}$ as $N\to\infty$, showing the tightness of \eqref
{eq:tight_seq} is equivalent to showing the tightness of the sequence
$\{Z_N'(t), t\in[-T,T]\}_{N\in\N}$, where $Z_N'$ is a process
defined by
%
\begin{equation}
Z_N'(t)=\frac{Z_N(t)-N}{N^{1/2}\mathrm{e}^{\psi(2)s_N/2}}.
\end{equation}

By a standard tightness criterion in the Skorokhod space given in \cite
{billingsleybook}, page~128, it suffices to show that there are $p>1$
and $C>0$ such that for all sufficiently large $N\in\N$ and all
$t_1,t_2,t_3\in[-T,T]$ with $t_1<t_2<t_3$,
%
\begin{equation}\label{eq:tighs_Need}
\E [|Z_N'(t_2)-Z_N'(t_1)|^{p}
|Z_N'(t_3)-Z_N'(t_2)|^{p}]\leq C|t_3-t_1|^{p}.
\end{equation}
It will be convenient to define random variables $X_1,\ldots ,X_N$ and
$Y_1,\ldots ,Y_N$ (which depend on $N,t_1,t_2,t_3$) by
\[
X_{i}=\mathrm{e}^{\xi_i(s_N+t_2)}-\mathrm{e}^{\xi_i(s_N+t_1)}, \qquad
Y_{i}=\mathrm{e}^{\xi_i(s_N+t_3)}-\mathrm{e}^{\xi_i(s_N+t_2)}.
\]
Then, we may rewrite \eqref{eq:tighs_Need} as follows:
%
\begin{equation}\label{eq:quadr_form_moment}
\E \Biggl|\sum_{i=1}^N \sum_{j=1}^N X_iY_j \Biggr|^{p}\leq C
N^p\mathrm{e}^{p\psi(2)s_N} |t_3-t_1|^p.
\end{equation}

First of all, we would like to treat the terms of the form $X_iY_i$ on
the left-hand side of \eqref{eq:quadr_form_moment} separately.
Applying Jensen's inequality $|\sum_{i=1}^k x_i|^p\leq k^{p-1}\sum
_{i=1}^k|x_i|^p$, $x_i\in\R$, we obtain
%
\begin{eqnarray}\label{eq:tight_jensen}
\E \Biggl|\sum_{i=1}^N \sum_{j=1}^N X_iY_j  \Biggr|^{p}
&=&\E \biggl|\sum_{1\leq i<j\leq N}X_iY_j+\sum_{1\leq j<i\leq
N}X_iY_j+\sum_{i=1}^N X_iY_i \biggr|^p \nonumber
\\[-8pt]
\\[-8pt]
&\leq&2\cdot3^{p-1} \E \biggl|\sum_{1\leq i<j\leq N}X_iY_j \biggr|^p+
3^{p-1}\E \Biggl|\sum_{i=1}^N X_iY_i \Biggr|^p.
\nonumber
\end{eqnarray}
In the rest of the proof we estimate the terms on the right-hand side.
We start by showing that
%
\begin{equation}\label{eq:tight_term1}
\E \Biggl|\sum_{i=1}^N X_iY_i \Biggr|^p\leq CN^p \mathrm{e}^{p\psi(2)s_N}|t_3-t_1|^p.
\end{equation}
By an inequality of Rosenthal \cite{rosenthal70}, Lemma 1 (or
see \cite{ibragimovsharakhmetov02}), 
%
\begin{equation}
\E \Biggl|\sum_{i=1}^N X_iY_i \Biggr|^p
\leq C\max \Biggl\{\sum_{i=1}^N\E|X_iY_i|^p,  \Biggl(\sum_{i=1}^N\E
|X_iY_i| \Biggr)^p \Biggr\}.
\end{equation}
Thus, to establish \eqref{eq:tight_term1}, it suffices to show that
%
\begin{eqnarray}\label
{eq:tight_term1a}
\E|X_iY_i|^p&\leq& CN^{p-1}\mathrm{e}^{p\psi(2)s_N}|t_3-t_1|^p,\\\label{eq:tight_term1b}
\E|X_iY_i|&\leq& C \mathrm{e}^{\psi(2)s_N}|t_3-t_1|.
\end{eqnarray}
Since $\xi$ is a process with stationary and independent increments,
we have
%
\begin{eqnarray}\label{eq:tight_wspom1}
\E|X_iY_i|^p
&=&
\E\bigl|\bigl(\mathrm{e}^{\xi(s_N+t_2)}-\mathrm{e}^{\xi(s_N+t_1)}\bigr)\bigl(\mathrm{e}^{\xi(s_N+t_3)}-\mathrm{e}^{\xi
(s_N+t_2)}\bigr)\bigr|^p \nonumber
\\[-8pt]
\\[-8pt]
&=&
\E\bigl[\mathrm{e}^{2p\xi(s_N+t_1)}\bigr] \cdot\E\bigl|\mathrm{e}^{\xi(t_3-t_2)}-1\bigr|^p\cdot\E
\bigl|\mathrm{e}^{2\xi(t_2-t_1)}-\mathrm{e}^{\xi(t_2-t_1)}\bigr|^p.
\nonumber
\end{eqnarray}
The first factor on the right-hand side of \eqref{eq:tight_wspom1}
equals $\mathrm{e}^{\psi(2p)(s_N+t_1)}$. Applying Lemma \ref{lem:small_incr}
to the last two factors on the right-hand side of \eqref
{eq:tight_wspom1}, we get
\[
\E|X_iY_i|^p\leq C \mathrm{e}^{\psi(2p)s_N}|t_3-t_1|^p.
\]
To complete the proof of \eqref{eq:tight_term1a}, we need to show that
for some $p>1$,
%
\begin{equation}\label{eq:tight_term1a_need}
\mathrm{e}^{(\psi(2p)-p\psi(2))s_N} \leq N^{p-1}.
\end{equation}
This is done as follows. Write for a moment $p=1+\delta$, where
$\delta>0$. By Assumption \eqref{eq:asymps_N1}, there is $\eps>0$
such that for sufficiently large $N$ we have $N^{p-1}>\mathrm{e}^{(\lambda
_2+\eps)\delta s_N}$. On the other hand, by Taylor's expansion,
\[
\psi(2p)-p\psi(2)=\delta\bigl(2\psi'(2)-\psi(2)\bigr)+\mathrm{o}(\delta)=\lambda
_2\delta+\mathrm{o}(\delta),  \qquad     \delta\to0,
\]
which is smaller than $(\lambda_2+\eps)\delta$ if $\delta$ is
sufficiently small.
Taking $\delta$ small enough, we~ob\-tain~\eqref{eq:tight_term1a_need}.
This completes the proof of \eqref{eq:tight_term1a}.

Let us prove \eqref{eq:tight_term1b}. Arguing as in \eqref
{eq:tight_wspom1}, we obtain
%
\begin{equation}\label{eq:tight_wspom1a}
\E|X_iY_i|
=
\E\bigl[\mathrm{e}^{2\xi(s_N+t_1)}\bigr] \cdot\E\bigl|\mathrm{e}^{\xi(t_3-t_2)}-1\bigr|\cdot\E\bigl|\mathrm{e}^{2\xi
(t_2-t_1)}-\mathrm{e}^{\xi(t_2-t_1)}\bigr|.
\end{equation}
The first factor on the right-hand side of \eqref{eq:tight_wspom1a}
equals $\mathrm{e}^{\psi(2)(s_N+t_1)}$. An application of Lemma \ref
{lem:small_incr} to the last two factors on the right-hand side
of \eqref{eq:tight_wspom1a} yields \eqref{eq:tight_term1b}.

We will now estimate the first term on the right-hand side of \eqref
{eq:tight_jensen}. We will show that
%
\begin{equation}\label{eq:tight_term2}
\E \biggl|\sum_{1\leq i<j\leq N}X_iY_j \biggr|^p\leq CN^p \mathrm{e}^{p\psi
(2)s_N}|t_3-t_1|^p.
\end{equation}
For $k=1,\ldots ,N$, denote by $\mathcal F_k$ the $\sigma$-algebra
generated by the random variables $X_1,\ldots ,X_k$ and $Y_1,\ldots
,Y_k$. Let $S_1=0$ and
%
\begin{equation}
S_k=\sum_{1\leq i<j\leq k} X_iY_j,   \qquad    k=2,\ldots ,N.
\end{equation}
We introduce also the sequence of differences $\Delta_1=0$ and
%
\begin{equation}\label{eq:tight_delta}
\Delta_k=S_k-S_{k-1}=Y_k(X_1+\cdots  +X_{k-1}),   \qquad   k=2,\ldots ,N.
\end{equation}
We claim that the sequence $\{S_k\}_{k=1}^N$ is a martingale with
respect to the filtration $\{\mathcal F_k\}_{k=1}^N$. Indeed, the
random variable $S_k$ is by definition $\mathcal F_k$-measurable, and
we have
\[
\E[S_k|\mathcal F_{k-1}]=S_{k-1}+\E[\Delta_k |\mathcal
F_{k-1}]=S_{k-1}+(X_1+\cdots  +X_{k-1})\E Y_k=S_{k-1},
\]
where the last equality follows from \eqref{eq:tight_simpl_ass}.
Having shown that $\{S_k\}_{k=1}^N$ is a martingale, we apply
Burkholder's inequality to obtain that for some constant $C=C(p)$,
%
\begin{equation}\label{eq:tight_burk}
\E|S_N|^p\leq C\E \Biggl(\sum_{i=1}^N \Delta_i^2 \Biggr)^{p/2}.
\end{equation}
The function $x\to x^{p/2}$, $x>0$, is concave since we choose $p$ to
be close to $1$. By Jensen's inequality applied to the right-hand side
of \eqref{eq:tight_burk},
%
\begin{equation}\label{eq:tight_jensen1}
\E|S_N|^p\leq C  \Biggl(\sum_{i=1}^N \E\Delta_i^2 \Biggr)^{p/2}.
\end{equation}
The random variables $Y_k$ and $X_1+\cdots  +X_{k-1}$ are independent,
and $\E X_k=0$, $k=1,\ldots ,N$, by \eqref{eq:tight_simpl_ass}. Hence,
by \eqref{eq:tight_delta}, $\E\Delta_k^2=(k-1)\E Y_1^2 \E X_1^2$. It
follows from \eqref{eq:tight_jensen1} that
%
\begin{equation}\label{eq:tight_wspom2}
\E|S_N|^p
\leq C  (N^2\E Y_1^2 \E X_1^2  )^{p/2}.
\end{equation}
We have, by Lemma \ref{lem:small_incr},
\[
\E X_1^2=\E\bigl[\mathrm{e}^{2\xi(s_N+t_1)}\bigr]\cdot\E\bigl(\mathrm{e}^{\xi(t_2-t_1)}-1\bigr)^2\leq C
\mathrm{e}^{\psi(2)s_N} (t_2-t_1).
\]
Similarly, $\E Y_1^2\leq C \mathrm{e}^{\psi(2)s_N} (t_3-t_2)$. Inserting this
into \eqref{eq:tight_wspom2}, we obtain
\[
\E|S_N|^p \leq C N^p \mathrm{e}^{p\psi(2)s_N}|t_3-t_1|^p.
\]
This proves \eqref{eq:tight_term2} and completes the proof of
tightness in Theorem \ref{theo:normal}.

\section{\texorpdfstring{Proof of Theorem \protect\ref{theo:crit}}{Proof of Theorem 1.3}}
\label{sec:proof_crit}
Let $W_N$ be a positive-valued stochastic process defined as in \eqref
{eq:def_bn}, that is,
%
\begin{equation}\label{eq:def_bn_crit}
W_N(t)=N^{-1/2}\mathrm{e}^{\xi(s_N+t)-\bfrac{\psi(2)}{2}(s_N+t)}.
\end{equation}
Fix $t_1\leq\cdots  \leq t_d$ and let $\mathbf{W}_{1,N},\ldots
,\mathbf{W}_{N,N}$ be independent copies of the $d$-dimensional random
vector $\mathbf{W}_{N}=(W_{N}(t_1),\ldots , W_N(t_d))$. Our aim is to
show that we have the following weak convergence of random vectors:
%
\begin{equation}\label{eq:theo_crit_eq_fd}
\sum_{i=1}^{N} (\mathbf{W}_{i,N}- \E\mathbf{W}_{i,N} )
\toweak
\bigl(\sqrt{\Phi(\vartheta)}\XXX(t_k)\bigr)_{k=1}^d,  \qquad     N\to\infty.
\end{equation}
In the one-dimensional case, the papers \cite
{benarousetal05,cranstonmolchanov05,kabluchkoprod} use the
classical summation theory of triangular arrays of random variables. We
will use a \textit{multidimensional} version of this theory
established in \cite{rvaceva62}; see \cite
{meerschaertschefflerbook} for a monograph treatment. According
to \cite{meerschaertschefflerbook}, Theorem~3.2.2 on page~53, we have
to verify that the following three conditions hold:
\begin{enumerate}[(3)]
\item[(1)] For every $\eps>0$,
%
\begin{equation}\label{eq:crit_cond1}
\lim_{N\to\infty} N\P[\|\mathbf{W}_{N}\|_{\infty}>\eps]=0.
\end{equation}
\item[(2)] For every $\eps>0$ and for every $\mathbf{v}=(v_1,\ldots
,v_d)\in\R^d$,
%
\begin{equation}\label{eq:crit_cond2}
\lim_{N\to\infty} N\Var \bigl[\langle\mathbf W_{N}, \mathbf
{v}\rangle1_{\|\mathbf{W}_{N}\|_{\infty}\leq\eps} \bigr]=
\Phi(\vartheta) \sum_{k,l=1}^d \mathrm{e}^{(\psi(1)-\afrac{\psi
(2)}{2})|t_l-t_k|}v_kv_l.
\end{equation}
\item[(3)] For every $\eps>0$,
%
\begin{equation}\label{eq:crit_cond3}
\lim_{N\to\infty} N \E\bigl[\mathbf{W}_{N} 1_{\|\mathbf{W}_{N}\|
_{\infty}> \eps}\bigr]=0.
\end{equation}
\end{enumerate}
Here, $\Phi$ is the standard normal distribution function and $\|\cdot
\|_{\infty}$ denotes the maximum norm on $\R^d$.

\subsection{\texorpdfstring{Proof of \protect\eqref{eq:crit_cond1} and
\protect\eqref{eq:crit_cond3}}{Proof of (66) and (68)}}

Let us first show that for every $t\in\R$ and every $\eps>0$, we have
%
\begin{equation}\label{eq:crit_wspom1}
\lim_{N\to\infty} N\E\bigl[W_N(t) 1_{W_{N}(t)> \eps}\bigr]=0.
\end{equation}
With $x_N=s_N+t$ and $b_N=\frac{1}{2}(\log N+\psi(2)x_N)+\log\eps$,
we may write
%
\begin{equation}\label{eq:crit_122}
N\E\bigl[W_N(t) 1_{W_{N}(t)> \eps}\bigr]
=N^{1/2}\mathrm{e}^{-\bfrac{\psi(2)}{2}x_N}\E\bigl[\mathrm{e}^{\xi(x_N)} 1_{\xi(x_N)>b_N}\bigr].
\end{equation}
Noting that by the critical growth condition \eqref{eq:asymps_Ncrit},
$\lim_{N\to\infty} b_N/x_N=\psi'(2)$ and applying part \hyperlink{p:2_trunc}{2} of Proposition \ref{prop:trunc_clt} with $\kappa=1$ to
the right-hand side of \eqref{eq:crit_122}, we obtain
%
\begin{eqnarray}\label{eq:crit_123}
N\E\bigl[W_N(t) 1_{W_{N}(t)> \eps}\bigr]
&\leq& C N^{1/2}\mathrm{e}^{-\bfrac{\psi(2)}{2}x_N} \mathrm{e}^{b_N}x_N^{-1/2}\mathrm{e}^{-I
(b_N/x_N )x_N}\nonumber
\\[-8pt]
\\[-8pt]
&\leq& C N x_N^{-1/2}\mathrm{e}^{-I (b_N/x_N )x_N}.
\nonumber
\end{eqnarray}
Using the convexity of the function $I$, as well as the fact that
$I(\psi'(2))=\lambda_2$ (see \eqref{eq:def_lambda}) and $I'(\psi
'(2))=2$ (see Lemma \ref{lem:I_prime_psi_prime}), we obtain
%
\begin{eqnarray}\label{eq:crit_124}
I \biggl(\frac{b_N}{x_N} \biggr)
&=&
I \biggl(\psi'(2)+\frac12  \biggl(\frac{\log N+2\log\eps
}{x_N}-\lambda_2 \biggr) \biggr)\notag\\
&\geq& I(\psi'(2))+I'(\psi'(2))\cdot\frac12  \biggl(\frac{\log
N+2\log\eps}{x_N}-\lambda_2 \biggr)\\
&=&
\frac{\log N+2\log\eps}{x_N}.\notag
\end{eqnarray}
It follows from \eqref{eq:crit_123} and \eqref{eq:crit_124} that
\[
N\E\bigl[W_N(t) 1_{W_{N}(t)> \eps}\bigr]\leq C N x_N^{-1/2}\mathrm{e}^{-\log N-2\log
\eps}\to0,  \qquad    N\to\infty.
\]
This proves \eqref{eq:crit_wspom1}. To prove \eqref{eq:crit_cond1},
note that
\[
N \P[\|\mathbf{W}_N\|_{\infty}>\eps]\leq N \sum_{k=1}^d \P
[W_N(t_k)>\eps]\leq \eps^{-1}N \sum_{k=1}^d \E
\bigl[W_N(t_k)1_{W_N(t_k)>\eps}\bigr].
\]
By \eqref{eq:crit_wspom1}, the right-hand side converges to $0$ as
$N\to\infty$. This proves \eqref{eq:crit_cond1}.

We proceed to the proof of \eqref{eq:crit_cond3}. Let $\Aa_{N,m}$,
$m=1,\ldots ,d$, be the random event $\{W_N(t_m)\geq W_N(t_l),
l=1,\ldots , d\}$. Then, for every $k=1,\ldots ,d$, we have
\begin{eqnarray*}
\E\bigl[W_{N}(t_k) 1_{\|\mathbf{W}_{N}\|_{\infty}> \eps}\bigr]
&\leq&\sum_{m=1}^d \E\bigl[W_{N}(t_k) 1_{\|\mathbf{W}_{N}\|_{\infty}>
\eps}1_{\Aa_{N,m}}\bigr]\\
&\leq&
\sum_{m=1}^d \E\bigl[W_N(t_m) 1_{W_{N}(t_m)>\eps}\bigr].
\end{eqnarray*}
An application of \eqref{eq:crit_wspom1} to the right-hand side
yields \eqref{eq:crit_cond3}.

\subsection{\texorpdfstring{Proof of \protect\eqref{eq:crit_cond2}}{Proof of (67)}}
\label{sec:proof_crit2}
It suffices to show that for every $1\leq k\leq l\leq d$ and every
$\eps>0$,
%
\begin{equation}\label{eq:crit_cond2b}
\lim_{N\to\infty} N \E\bigl[W_N(t_k)W_N(t_l) 1_{\|\mathbf{W}_N\|
_{\infty}\leq\eps}\bigr]=\Phi(\vartheta) \mathrm{e}^{(\psi(1)-\afrac{\psi
(2)}{2})(t_l-t_k)}.
\end{equation}
Let us start by computing a closely related limit. We will show that
%
\begin{equation}\label{eq:crit1}
\lim_{N\to\infty}N\E\bigl[W_{N}(t_k)W_{N}(t_l)1_{W_N(t_1)\leq\eps}\bigr]
=
\Phi(\vartheta) \mathrm{e}^{(\psi(1)-\afrac{\psi(2)}{2})(t_l-t_k)}.
\end{equation}
It follows from \eqref{eq:def_bn_crit} that
%
\begin{equation}\label{eq:crit_125}
\E\bigl[W_{N}(t_k)W_{N}(t_l)1_{W_N(t_1)\leq\eps}\bigr]
=\frac{\E[\mathrm{e}^{\xi(s_N+t_k)+\xi(s_N+t_l)}
1_{W_N(t_1)\leq\eps}]}{N\mathrm{e}^{\psi(2)s_N}\mathrm{e}^{\bfrac{\psi(2)}{2}(t_k+t_l)}}.
\end{equation}
Using the fact that $\xi$ is a L\'evy process, we obtain
%
\begin{eqnarray}\label{eq:crit_126}
&&\E\bigl[\mathrm{e}^{\xi(s_N+t_k)+\xi(s_N+t_l)}1_{W_N(t_1)\leq\eps
}\bigr]\notag\\
&& \quad =
\E\bigl[\mathrm{e}^{2\xi(s_N+t_1)}1_{W_N(t_1)\leq\eps}\bigr]\cdot
\E \mathrm{e}^{\xi(s_N+t_k)+\xi(s_N+t_l)-2\xi(s_N+t_1)}\nonumber
\\[-8pt]
\\[-8pt]
&& \quad =
\E\bigl[\mathrm{e}^{2\xi(s_N+t_1)}1_{W_N(t_1)\leq\eps}\bigr]\cdot\E \mathrm{e}^{\xi
(t_k-t_1)+\xi(t_l-t_1)}\notag\\
&& \quad =
\E\bigl[\mathrm{e}^{2\xi(x_N)}1_{\xi(x_N)\leq b_N}\bigr]\cdot\E \mathrm{e}^{\xi(t_k-t_1)+\xi
(t_l-t_1)},\notag
\end{eqnarray}
where we have used the notation
%
\begin{equation}\label{eq:pr_crit_xN}
x_N=s_N+t_1,  \qquad    b_N=\tfrac12 \bigl(\log N+\psi(2)x_N\bigr)+\log\eps.
\end{equation}
The critical growth condition \eqref{eq:asymps_Ncrit} implies that
%
\begin{equation}
b_N=\psi'(2)x_N+\vartheta\sqrt{\psi''(2)x_N}+
\mathrm{o}\bigl(\sqrt{x_N}\bigr),\qquad
N\to\infty.
\end{equation}
Applying part \hyperlink{p:1_trunc}{1} of Proposition \ref{prop:trunc_clt}
with $\kappa=2$, we obtain
%
\begin{equation}\label{eq:crit_127}
\E\bigl[\mathrm{e}^{2\xi(x_N)}1_{\xi(x_N)\leq b_N}\bigr] \sim\Phi(\vartheta) \mathrm{e}^{\psi
(2) (s_N+t_1)},  \qquad    N\to\infty.
\end{equation}
Recalling that $\xi$ is a L\'evy process and taking into account that
$t_k\leq t_l$, we obtain
%
\begin{equation}\label{eq:crit_127a}
\E \mathrm{e}^{\xi(t_k-t_1)+\xi(t_l-t_1)}=\mathrm{e}^{\psi(2)(t_k-t_1)}\mathrm{e}^{\psi(1)(t_l-t_k)}.
\end{equation}
Bringing equations \eqref{eq:crit_125}, \eqref{eq:crit_126}, \eqref
{eq:crit_127} and \eqref{eq:crit_127a} together, we obtain \eqref{eq:crit1}.
Trivially, it follows from \eqref{eq:crit1} that
%
\begin{equation}
\limsup_{N\to\infty}N\E\bigl[W_{N}(t_k)W_{N}(t_l)1_{\|\mathbf{W}_N\|
_{\infty}\leq\eps}\bigr]
\leq
\Phi(\vartheta) \mathrm{e}^{(\psi(1)-\afrac{\psi(2)}{2})(t_l-t_k)}.
\end{equation}
We are going to prove the converse inequality:
%
\begin{equation}
\liminf_{N\to\infty}N\E\bigl[W_{N}(t_k)W_{N}(t_l)1_{\|\mathbf{W}_N\|
_{\infty}\leq\eps}\bigr]
\geq
\Phi(\vartheta) \mathrm{e}^{(\psi(1)-\afrac{\psi(2)}{2})(t_l-t_k)}.
\end{equation}
Note that for every (small) $\eta>0$, the following inclusion of
random events holds:
\[
\{\|\mathbf{W}_N\|_{\infty}\leq\eps\}
\supset
\{W_N(t_1)\leq\eta\eps\} \Big  \backslash\bigcup_{m=1}^d \Aa_{N,m},
\]
where $\Aa_{N,m}$ is the random event $\{\xi(s_N+t_m)-\xi
(s_N+t_1)>-\log\eta\}$. Thus,
\begin{eqnarray*}
&& \E\bigl[W_{N}(t_k)W_{N}(t_l)1_{\|\mathbf{W}_N\|_{\infty}\leq
\eps}\bigr] \\
&& \quad \geq
\E\bigl[W_{N}(t_k)W_{N}(t_l)1_{W_N(t_1)\leq\eta\eps}\bigr]-
\sum_{m=1}^d \E[W_{N}(t_k)W_{N}(t_l)1_{\Aa_{N,m}}].
\end{eqnarray*}
Since the asymptotic behavior of the first term on the right-hand side
was computed in \eqref{eq:crit1}, we need to show that for every
$m=1,\ldots ,d$, and every $1\leq k\leq l\leq d$,
%
\begin{equation}\label{eq:crit_128}
\lim_{\eta\downarrow0}\limsup_{N\to\infty} N\E
[W_{N}(t_k)W_{N}(t_l)1_{\Aa_{N,m}}]=0.
\end{equation}
By \eqref{eq:def_bn_crit}, we have
%
\begin{eqnarray}\label
{eq:crit_129}
&& \E[W_{N}(t_k)W_{N}(t_l)1_{\Aa_{N,m}}] \notag \\
&& \quad \leq
CN^{-1} \mathrm{e}^{-\psi(2)s_N} \E\bigl[\mathrm{e}^{\xi(s_N+t_k)+\xi(s_N+t_l)}1_{\Aa
_{N,m}}\bigr]\nonumber
\\[-8pt]
\\[-8pt]
&& \quad =
CN^{-1} \mathrm{e}^{-\psi(2)s_N} \E\bigl[\mathrm{e}^{2\xi(s_N+t_1)}\mathrm{e}^{\xi(s_N+t_k)+\xi
(s_N+t_l)-2\xi(s_N+t_1)}1_{\Aa_{N,m}}\bigr]\notag\\
&& \quad \leq CN^{-1}\E\bigl[\mathrm{e}^{\xi(t_k-t_1)+\xi(t_l-t_1)}1_{\xi(t_m-t_1)>-\log
\eta}\bigr].\notag
\end{eqnarray}
Note that by \eqref{eq:def_varphi}, $\E \mathrm{e}^{\xi(t_k-t_1)+\xi
(t_l-t_1)}<\infty$. Hence, by the dominated convergence theorem,
%
\begin{equation}\label{eq:crit_130}
\lim_{\eta\downarrow0} \E\bigl[\mathrm{e}^{\xi(t_k-t_1)+\xi(t_l-t_1)}1_{\xi
(t_m-t_1)>-\log\eta}\bigr]=0.
\end{equation}
To complete the proof of \eqref{eq:crit_128}, combine \eqref
{eq:crit_129} and \eqref{eq:crit_130}.

\section{\texorpdfstring{Proof of Theorem \protect\ref{theo:stable}}{Proof of Theorem 1.4}}
\label{sec:proof_stable}
\subsection{Notation and preliminaries}
We will concentrate on proving the convergence in the Skorokhod space
$D[0,T]$. For the proof of the two-sided convergence on $D[-T,T]$ we
refer to \cite{kabluchkoprep}.

We start by introducing some notation. Let $W_{1,N},\ldots ,W_{N,N}$ be
independent copies of a positive-valued random process $\{
W_{N}(t),t\geq0\}$ defined by
%
\begin{equation}\label{eq:def_wn_1}
W_{N}(t)=\mathrm{e}^{\xi(s_N+t)-b_N(t)},
\end{equation}
where $b_N(t)$ is given by
%
\begin{equation}\label{eq:def_bn_stable}
b_N(t)= \log B_N(t)=\frac{\psi(\alpha)}{\alpha}t+s_NI^{-1}
\biggl(\frac{\log N-\log( \alpha\sqrt{2\uppi\psi''(\alpha) s_N})
}{s_N} \biggr).
\end{equation}
Define a process $Y_N$ by
%
\begin{equation}\label{eq:def_YN}
Y_{N}(t)=\frac{Z_N(t)-A_N(t)}{B_N(t)}=
\cases{\displaystyle
\sum_{i=1}^N W_{i,N}(t),
&\quad$0<\alpha<1$,\cr\displaystyle
\sum_{i=1}^N W_{i,N}(t)-N\E\bigl[W_{N}(t)1_{W_{N}(0)\leq1}\bigr],
&\quad$\alpha=1$,\cr\displaystyle
\sum_{i=1}^N W_{i,N}(t) -N\E W_{N}(t),
&\quad$1<\alpha<2$.
}
\end{equation}
Our aim is to show that we have the following weak convergence of
stochastic processes on the Skorokhod space $D[0,T]$:
%
\begin{equation}\label{eq:pr_stab_YN_YYY}
Y_N(\cdot)\toweak\YYY_{\alpha;\xi}(\cdot),  \qquad   N\to\infty.
\end{equation}

We will use an approach based on considering the extremal order
statistics. This method goes back to LePage \textit{et al.} \cite{lepageetal81} and was
used in the context of the random energy model by Bovier \textit{et al.} \cite{bovieretal02} (note that the papers \cite
{benarousetal05,cranstonmolchanov05,kabluchkoprod} use a different
method). To describe the method of our proof of \eqref
{eq:pr_stab_YN_YYY}, let us consider the case $\alpha\in(0,1)$ only.
The first step is to prove that the upper order statistics of the
sequence $W_{1,N}(0),\ldots ,W_{N,N}(0)$ can be approximated, as $N\to
\infty$, by the Poisson process $\{U_i,i\in\N\}$ defined as in
Section~\ref{subsec:def_proc_Y_axi}.
In the second step we write, for $t\geq0$,
%
\begin{equation}\label{eq:def_eta_iN}
\sum_{i=1}^N W_{i,N}(t)
=\sum_{i=1}^N W_{i,N}(0) \mathrm{e}^{\eta_{i,N}(t)},
\end{equation}
where $\{\eta_{i,N}(t), t\geq0\}$, $i=1,\ldots ,N$, are processes
defined by
%
\begin{equation}\label{eq:def_eta_iN1}
\eta_{i,N}(t)=\xi_i(s_N+t)-\xi_i(s_N)-\frac{\psi(\alpha)}{\alpha}t.
\end{equation}
Note that the processes $\eta_{1,N},\ldots ,\eta_{N,N}$ are
independent of each other, independent of $W_{1,N}(0), \ldots
,W_{N,N}(0)$, and have the same law as the process $\eta$ defined by
$\eta(t)=\xi(t)-\frac{\psi(\alpha)}{\alpha}t$. Bringing
everything together, we may write
%
\begin{equation}\label{eq:pr_stab_heuristic}
\sum_{i=1}^N W_{i,N}(t)\to
\sum_{i=1}^{\infty}U_i \mathrm{e}^{\xi_i(t)-\bfrac{\psi(\alpha)}{\alpha
}t}=\YYY_{\alpha;\xi}(t),  \qquad    N\to\infty.
\end{equation}
%
The rest of the section is devoted to the justification of the above argument.

\subsection{Asymptotics for truncated moments}
The following corollary of Proposition \ref{prop:trunc_clt} will play
a crucial role in the sequel.
\begin{proposition}\label{lem:aux_int}
Let the assumptions of Theorem \ref{theo:stable} be satisfied. Let
$W_N$ be a process defined by \eqref{eq:def_wn_1}. The following three
statements hold true.
\begin{enumerate}[(3)]
\item[(1)]\hypertarget{p:2_lem6} Let $0\leq\kappa<\alpha$. Then, for every
$\tau>0$,
%
\begin{equation}\label{eq:lem_aux_3}
\lim_{N\to\infty} N \E \bigl[W_N^{\kappa}(0)1_{W_N(0)> \tau
} \bigr]=\frac{\alpha}{\alpha-\kappa} \tau^{\kappa-\alpha}.
\end{equation}
\item[(2)]\hypertarget{p:1_lem6} Let $\kappa>\alpha$. Then, for every $\tau>0$,
%
\begin{equation}\label{eq:lem_aux_2}
\lim_{N\to\infty} N \E \bigl[W_N^{\kappa}(0)1_{W_N(0)\leq\tau
} \bigr]=\frac{\alpha}{\kappa-\alpha} \tau^{\kappa-\alpha}.
\end{equation}
\item[(3)]\hypertarget{p:3_lem6} Let $\kappa=\alpha$. Then, for every $0<\tau
_1\leq\tau_2$,
%
\begin{equation}\label{eq:lem_aux_4}
\lim_{N\to\infty} N \E \bigl[W_N^{\kappa}(0)1_{W_N(0)\in (\tau
_1,\tau_2)} \bigr]=\kappa(\log\tau_2-\log\tau_1).
\end{equation}
\end{enumerate}
\end{proposition}
\begin{pf}
We prove part \hyperlink{p:2_lem6}{1} of the proposition. Recall from \eqref
{eq:def_bn_stable} that
%
\begin{equation}
b_N(0)=s_NI^{-1}(c_N),  \qquad \mbox{where } c_N=\frac{\log N-\log( \alpha
\sqrt{2\uppi\psi''(\alpha) s_N}) }{s_N}.
\end{equation}
We have $\lim_{N\to\infty}I^{-1}(c_N)=\psi'(\alpha)$ by the fast
growth condition \eqref{eq:asymps_N}.
By part \hyperlink{p:2_trunc}{2} of Proposition~\ref{prop:trunc_clt}, we have
as $N\to\infty$,
%
\begin{eqnarray}\label{eq:pr_stab_moment1}
\E \bigl[W_N^{\kappa}(0)1_{W_N(0)> \tau} \bigr]
&=&\mathrm{e}^{-\kappa b_N(0)}\E\bigl[\mathrm{e}^{\kappa\xi(s_N)}1_{\xi(s_N)>b_N(0)+\log
\tau}\bigr] \nonumber
\\[-8pt]
\\[-8pt]
&\sim&
\frac{\tau^{\kappa}}{(\alpha-\kappa) \sqrt{2\uppi\psi''(\alpha
)s_N}} \mathrm{e}^{-I((b_N(0)+\log\tau)/s_N)s_N}.\notag
\end{eqnarray}
To compute the asymptotic behavior of the right-hand side of \eqref
{eq:pr_stab_moment1}, we will prove that
%
\begin{equation}\label{eq:pr_stab_sN_I}
s_NI \biggl(\frac{b_N(0)+\log\tau}{s_N} \biggr)=s_Nc_N+\alpha\log
\tau+\mathrm{o}(1),   \qquad   N\to\infty.
\end{equation}
We have $\lim_{N\to\infty}I^{-1}(c_N)=\psi'(\alpha)$, hence $\lim
_{N\to\infty}I'(I^{-1}(c_N))=\alpha$ by Lemma \ref
{lem:I_prime_psi_prime}. Using Taylor's expansion of $I$ around the
point $I^{-1}(c_N)$, we obtain
\[
I \biggl(\frac{b_N(0)+\log\tau}{s_N} \biggr)
=I \biggl(I^{-1} (c_N )+\frac{\log\tau}{s_N} \biggr)
=c_N+\frac{\alpha\log\tau+\mathrm{o}(1)}{s_N},  \qquad     N\to\infty.
\]
This proves \eqref{eq:pr_stab_sN_I}. Inserting \eqref
{eq:pr_stab_sN_I} into \eqref{eq:pr_stab_moment1}, we obtain part \hyperlink{p:2_lem6}{1}
 of the proposition. Part \hyperlink{p:1_lem6}{2} can be proved in a
similar way.

Let us prove part \hyperlink{p:3_lem6}{3}  of the proposition. We write
$F_N(\tau)=\P[W_N(0)\leq\tau]$ for the distribution function of
$W_N(0)$, and $\bar F_N(\tau)=1-F_N(\tau)$ for its tail. Taking
$\kappa=0$ in \eqref{eq:lem_aux_3}, we obtain
%
\begin{equation}\label{eq:pr_stab_trunc_tail}
\lim_{N\to\infty}N\bar F_N(\tau)=\tau^{-\alpha}.
\end{equation}
Note that this holds uniformly in $\tau\in(\tau_1,\tau_2)$, cf.\
Theorem \ref{theo:ld}. Trivially, we have
\[
N \E \bigl[W_N^{\kappa}(0)1_{W_N(0)\in (\tau_1,\tau_2)}
\bigr]=N\int_{\tau_1}^{\tau_2} w^{\kappa}\,\mathrm{d}F_N(w)=-N\int_{\tau
_1}^{\tau_2} w^{\kappa}\,\mathrm{d}\bar F_N(w).
\]
Integrating by parts, we obtain
\[
N \E \bigl[W_N^{\kappa}(0)1_{W_N(0)\in (\tau_1,\tau_2)}
\bigr]= -w^{\kappa}N\bar F_N(w) |_{\tau_1}^{\tau_2}+\kappa
\int_{\tau_1}^{\tau_2}w^{\kappa-1}N\bar F_N(w)\,\mathrm{d}w.
\]
Applying \eqref{eq:pr_stab_trunc_tail} to the right-hand side and
recalling that $\kappa=\alpha$, we obtain
\[
\lim_{N\to\infty} N \E \bigl[W_N^{\kappa}(0)1_{W_N(0)\in (\tau
_1,\tau_2)} \bigr]=\kappa\int_{\tau_1}^{\tau_2}w^{-1}\,\mathrm{d}w=\kappa
(\log\tau_2-\log\tau_1),
\]
which completes the proof of part \hyperlink{p:3_lem6}{3}.
\end{pf}


\subsection{Convergence of the upper order statistics}
For $\tau>0$, we define a process $\YYY_{\alpha;\xi}^{(\tau,\infty
)}$, which is a ``truncated version'' of the process $\YYY_{\alpha
;\xi}$, by
%
\begin{equation}\label{eq:def_Ytau}
\YYY_{\alpha;\xi}^{(\tau,\infty)}(t)=
\cases{\displaystyle
\mathop{\mathop{\sum}_{i\in\N}}_{U_i>\tau} U_i \mathrm{e}^{\xi_i(t)-\bfrac{\psi
(\alpha)}{\alpha}t},
&\quad$0<\alpha<1$,\cr\displaystyle
\mathop{\mathop{\sum}_{i\in\N}}_{U_i>\tau}U_i \mathrm{e}^{\xi_i(t)-\psi(1)t}-\log
\frac{1}{\tau},
&\quad$\alpha=1$,\cr\displaystyle
\mathop{\mathop{\sum}_{i\in\N}}_{ U_i>\tau} U_i \mathrm{e}^{\xi_i(t)-\bfrac{\psi
(\alpha)}{\alpha}t}
-\frac{\alpha \tau^{1-\alpha}}{\alpha-1} \mathrm{e}^{(\psi(1)-\bfrac{\psi
(\alpha)}{\alpha})t},
&\quad$1<\alpha<2$.
}
\end{equation}
Similarly, we define $Y_{N}^{(\tau,\infty)}$, a truncated version of
the process $Y_N$ given by \eqref{eq:def_YN}, by
%
\begin{equation}\label{eq:def_YNtau}
Y_{N}^{(\tau,\infty)}(t)=
\cases{\displaystyle
\mathop{\mathop{\sum}_{1\leq i\leq N}}_{W_{i,N}(0)>\tau} W_{i,N}(t),
&\quad$0<\alpha<1$,\cr\displaystyle
\mathop{\mathop{\sum}_{1\leq i\leq N}}_{W_{i,N}(0)>\tau} W_{i,N}(t)-N\E
\bigl[W_{N}(t)1_{W_{N}(0)\in(\tau,1)}\bigr],
&\quad$\alpha=1$,\cr\displaystyle
\mathop{\mathop{\sum}_{1\leq i\leq N}}_{W_{i,N}(0)>\tau} W_{i,N}(t) -N\E
\bigl[W_{N}(t)1_{W_{N}(0)>\tau}\bigr],
&\quad$1<\alpha<2$.
}
\end{equation}
The next lemma is the main result of this subsection.
\begin{lemma}\label{lem:est_upper_stat}
For every $\tau>0$, we have the following weak convergence of
stochastic processes on the Skorokhod space $D[0,T]$:
\[
Y^{(\tau,\infty)}_N(\cdot)
\toskor
\YYY^{(\tau,\infty)}_{\alpha;\xi}(\cdot),  \qquad    N\to\infty.
\]
\end{lemma}

First, we establish the convergence of regularizing terms in \eqref
{eq:def_YNtau} to those in \eqref{eq:def_Ytau}. If $\alpha\in(1,2)$,
then writing $W_{N}(t)=W_{N}(0)\mathrm{e}^{\eta_{N}(t)}$ with $\eta_N(t)=\xi
(s_N+t)-\xi(s_N)-\frac{\psi(\alpha)}{\alpha}t$ (see
equations \eqref{eq:def_eta_iN} and \eqref{eq:def_eta_iN1}) and
applying part \hyperlink{p:2_lem6}{1} of Proposition \ref{lem:aux_int}, we obtain
\begin{eqnarray*}
\lim_{N\to\infty} N\E\bigl[W_{N}(t)1_{W_{N}(0)>\tau}\bigr]
&=&
\mathrm{e}^{(\psi(1)-\afrac{\psi(\alpha)}{\alpha})t}\lim_{N\to\infty} N\E
\bigl[W_{N}(0)1_{W_{N}(0)>\tau}\bigr]\\
&=&
\frac{\alpha\tau^{1-\alpha}}{\alpha-1} \mathrm{e}^{(\psi(1)-\afrac{\psi
(\alpha)}{\alpha})t}.
\end{eqnarray*}
If $\alpha=1$, then part \hyperlink{p:3_lem6}{3} of Proposition \ref
{lem:aux_int} yields
\[
\lim_{N\to\infty} N\E\bigl[W_{N}(t)1_{W_{N}(0)\in(\tau,1)}\bigr]=\lim
_{N\to\infty} N\E\bigl[W_{N}(0)1_{W_{N}(0)\in(\tau,1)}\bigr]=\log\frac
1{\tau}.
\]

Thus, in proving Lemma \ref{lem:est_upper_stat}, we may drop the
regularizing terms in \eqref{eq:def_Ytau} and \eqref{eq:def_YNtau}.
More precisely, we define stochastic processes $\tilde\YYY_{\alpha
;\xi}^{(\tau,\infty)}$ and $\tilde Y^{(\tau,\infty)}_{N}$ by
%
\begin{eqnarray}\label{eq:def_Y_tilde_tau}
\tilde\YYY_{\alpha;\xi}^{(\tau,\infty)}(t)
&=&
\mathop{\mathop{\sum}_{i\in\N}}_{U_i>\tau} U_i \mathrm{e}^{\xi_i(t)-\bfrac{\psi
(\alpha)}{\alpha}t},\\\label{eq:def_YN_tilde_tau}
\tilde Y^{(\tau,\infty)}_{N}(t)
&=&
\mathop{\mathop{\sum}_{1\leq i\leq N}}_{W_{i,N}(0)>\tau} W_{i,N}(t)
=
\mathop{\mathop{\sum}_{1\leq i\leq N}}_{W_{i,N}(0)>\tau} W_{i,N}(0)\mathrm{e}^{\eta
_{i,N}(t)};
\end{eqnarray}
see \eqref{eq:def_eta_iN} and \eqref{eq:def_eta_iN1} for the last
equality. With this notation, we may restate Lemma \ref
{lem:est_upper_stat} as follows.
\begin{lemma}\label{lem:est_upper_stat1}
For every $\tau>0$, we have the following weak convergence of
stochastic processes on the Skorokhod space $D[0,T]$:
%
\begin{equation}\label{eq:upper_stat_conv_proc}
\tilde Y^{(\tau,\infty)}_N(\cdot)
\toskor
\tilde\YYY^{(\tau,\infty)}_{\alpha;\xi}(\cdot),  \qquad    N\to
\infty.
\end{equation}
\end{lemma}

We start by considering the upper order statistics of the summands on
the right-hand side of \eqref{eq:def_YN_tilde_tau} at $t=0$. More
precisely, let $\{W_{i:N}(0)\}_{i=1}^{N}$ be the rearrangement of the
numbers $\{W_{i,N}(0)\}_{i=1}^N$ in the descending order, and set also
$W_{i:N}(0)=0$ for $i>N$. Let $\SSS$ be the space of all sequences
$w=(w_i)_{i=1}^{\infty}$ with $w_1\geq w_2\geq\cdots  \geq0$. Then,
$\SSS$ is a closed subset of $\R^{\infty}$ endowed with the product topology.
\begin{lemma}\label{lem:order_stat}
Let $\{U_i, i\in\N\}$ be the points of a Poisson process on
$(0,\infty)$ with intensity $\alpha u^{-(\alpha+1)}\,\mathrm{d}u$, arranged in
the descending order. Then, we have the following weak convergence of
random elements in $\SSS$:
%
\begin{equation}\label{eq:upper_stat_conv}
\{W_{i:N}(0)\}_{i=1}^{\infty}\toweak\{U_i\}_{i=1}^{\infty}, \qquad
N\to\infty.
\end{equation}
\end{lemma}
\begin{pf}
By part \hyperlink{p:2_lem6}{1} of Proposition \ref{lem:aux_int} with $\kappa
=0$, we have for every $u>0$,
%
\begin{equation}
\lim_{N\to\infty}N\P[W_{N}(0)>u]=u^{-\alpha}.
\end{equation}
%
To complete the proof, use \cite{resnickbook}, Proposition 3.21 on page~154.
\end{pf}
\begin{pf*}{Proof of Lemma \ref{lem:est_upper_stat1}}
Let $f\dvtx D[0,T]\to\R$ be a continuous bounded function. To pro\-ve~\eqref
{eq:upper_stat_conv_proc}, we need to verify that
%
\begin{equation}\label{eq:pst101}
\lim_{N\to\infty} \E f\bigl(\tilde Y^{(\tau,\infty)}_N\bigr)=\E f\bigl(\tilde
\YYY^{(\tau,\infty)}_{\alpha;\xi}\bigr).
\end{equation}
Let $\SSS_{\tau}\subset\SSS$ be the set of all sequences
$(w_i)_{i\in\N}\in\SSS$ with $\lim_{i\to\infty}w_i=0$ and such
that $w_i\neq\tau$ for all $i\in\N$. Define a function $\bar f\dvtx \SSS
_{\tau}\to\R$ by
\[
\bar f(w)=\E f \biggl(\mathop{\mathop{\sum}_{i\in\N}}_{w_i>\tau} w_i \mathrm{e}^{\xi
_i(\cdot)-\afrac{\psi(\alpha)}{\alpha}\cdot} \biggr), \qquad
w=(w_i)_{i\in\N}\in\SSS_{\tau}.
\]
Note that $\bar f$ is bounded and continuous on $\SSS_{\tau}$, and
$\SSS_{\tau}$ has full measure with respect to the law of
$(U_i)_{i=1}^{\infty}$. By Fubini's theorem,
%
\begin{equation}\label{eq:pst100}
\E f\bigl(\tilde Y^{(\tau,\infty)}_N\bigr)=\E\bar f((W_{i:N}(0))_{i=1}^{\infty
}), \qquad
\E f\bigl(\tilde\YYY^{(\tau,\infty)}_{\alpha;\xi}\bigr)=\E\bar
f((U_i)_{i=1}^{\infty}).
\end{equation}
It follows from Lemma \ref{lem:order_stat} and the properties of the
weak convergence that
%
\begin{equation}\label{eq:pst102}
\lim_{N\to\infty}\E\bar f((W_{i:N}(0))_{i=1}^{\infty})=\E\bar
f((U_i)_{i=1}^{\infty}).
\end{equation}
Putting \eqref{eq:pst100} and \eqref{eq:pst102} together, we
obtain \eqref{eq:pst101}. This completes the proof of the lemma.
\end{pf*}

\subsection{Estimating the lower order statistics}
In this section we estimate the difference between the processes $\YYY
_{\alpha;\xi}$ and $Y_N$ and their truncated versions $\YYY_{\alpha
;\xi}^{(\tau,\infty)}$ and $Y_N^{(\tau,\infty)}$. Define a process
$\YYY_{\alpha;\xi}^{(0,\tau)}$ by
%
\begin{equation}\label{eq:def_Y_zero_tau}
\YYY_{\alpha;\xi}^{(0,\tau)}(t)=\YYY_{\alpha;\xi}(t)-\YYY
_{\alpha;\xi}^{(\tau,\infty)}(t).
\end{equation}
\begin{lemma}\label{lem:est_lower_stat1}
For every $\eps>0$, we have
%
\begin{equation}
\lim_{\tau\downarrow0} \P\Bigl [\sup_{t\in[0,T]} \bigl|\YYY_{\alpha
;\xi}^{(0,\tau)}(t)\bigr|>\eps \Bigr]=0.
\end{equation}
\end{lemma}
\begin{pf}
The proof follows immediately from Proposition \ref{prop:uniform}.
\end{pf}

Next we define a process $Y^{(0,\tau)}_{N}$ representing the sum of
the lower order statistics in~\eqref{eq:def_YN} by $Y^{(0,\tau
)}_{N}(t)=Y_N(t)-Y_N^{(\tau,\infty)}(t)$. Equivalently,
%
\begin{equation}\label{eq:def_Y_lower}
Y^{(0,\tau)}_{N}(t)=
\cases{\displaystyle
\mathop{\mathop{\sum}_{1\leq i\leq N}}_{ W_{i,N}(0)\leq\tau} W_{i,N}(t),
&\quad$\alpha\in(0,1)$,\cr\displaystyle
\mathop{\mathop{\sum}_{1\leq i\leq N}}_{ W_{i,N}(0)\leq\tau} W_{i,N}(t)-N\E
\bigl[W_{N}(t)1_{W_{N}(0)\leq\tau}\bigr],
&\quad$\alpha\in[1,2)$.
}
\end{equation}
\begin{lemma}\label{lem:est_lower_stat2}
For every $\eps>0$, we have
%
\begin{equation}
\lim_{\tau\downarrow0}\limsup_{N\to\infty} \P \Bigl[\sup_{t\in
[0,T]} \bigl|Y_N^{(0,\tau)}(t)\bigr|>\eps \Bigr]=0.
\end{equation}
\end{lemma}

The proof will be carried out in the rest of the subsection. First we
consider the regularizing term in \eqref{eq:def_Y_lower}.
If $\alpha\in(0,1)$, then applying part \hyperlink{p:1_lem6}{2} of
Proposition \ref{lem:aux_int} with $\kappa=1$, we obtain
%
\begin{equation}\label{eq:pr_stab_reg_term1}
\lim_{\tau\downarrow0}\limsup_{N\to\infty}N\E
\bigl[W_{N}(t)1_{W_{N}(0)\leq\tau}\bigr]
=0.
\end{equation}
Define a process $\tilde Y^{(0,\tau)}_{N}$ coinciding with $Y^{(0,\tau
)}_{N}$ for $\alpha\in[1,2)$ and containing an additional term for
$\alpha\in(0,1)$ by
%
\begin{equation}\label{eq:def_Y_tilde_lower}
\tilde Y^{(0,\tau)}_{N}(t)=
\mathop{\mathop{\sum}_{1\leq i\leq N}}_{ W_{i,N}(0)\leq\tau} W_{i,N}(t)-N\E
\bigl[W_{N}(t)1_{W_{N}(0)\leq\tau}\bigr].
\end{equation}
In view of \eqref{eq:pr_stab_reg_term1}, we may restate Lemma \ref
{lem:est_lower_stat2} as follows.
\begin{lemma}
For every $\eps>0$, we have
%
\begin{equation}
\lim_{\tau\downarrow0}\limsup_{N\to\infty} \P \Bigl[\sup_{t\in
[0,T]} \bigl|\tilde Y_N^{(0,\tau)}(t)\bigr|>\eps \Bigr]=0.
\end{equation}
%
\end{lemma}
\begin{pf}
For a function $f\dvtx [0,T]\to\R$ we write $\|f\|_{\infty}=\sup_{t\in
[0,T]}|f(t)|$. We have
%
\begin{eqnarray}\label
{eq:pr_wspom333}
\tilde Y_N^{(0,\tau)}(t)
&=&
\sum_{i=1}^N \bigl(W_{i,N}(0)1_{W_{i,N}(0)\leq\tau}-\E
\bigl[W_{N}(0)1_{W_{N}(0)\leq\tau}\bigr]\bigr) \mathrm{e}^{\eta_{i,N}(t)}\nonumber
\\[-8pt]
\\[-8pt]
&&{}+
\E\bigl[W_{N}(0)1_{W_{N}(0)\leq\tau}\bigr] \sum_{i=1}^N \bigl(\mathrm{e}^{\eta
_{i,N}(t)}-\E \mathrm{e}^{\eta_{i,N}(t)}\bigr).\notag
\end{eqnarray}
It follows from \eqref{eq:pr_wspom333} that $\|\tilde Y_N^{(0,\tau)}\|
_{\infty}\leq M_{N,\tau}'+M_{N,\tau}''$, where $M_{N,\tau}'$ and
$M_{N,\tau}''$ are random variables defined by
\begin{eqnarray*}
M_{N,\tau}'&=&\sum_{i=1}^N \|\mathrm{e}^{\eta_{i,N}}\|_{\infty}
\bigl|W_{i,N}(0)1_{W_{i,N}(0)\leq\tau}-\E\bigl[W_{N}(0)1_{W_{N}(0)\leq\tau}\bigr]\bigr|
,\\
M_{N,\tau}''&=&\E\bigl[W_{N}(0)1_{W_{N}(0)\leq\tau}\bigr] \cdot  \Biggl\|\sum
_{i=1}^N (\mathrm{e}^{\eta_{i,N}}-\E \mathrm{e}^{\eta_{i,N}}) \Biggr\|_{\infty}.
\end{eqnarray*}
Thus, to prove the lemma, it suffices to show that
%
\begin{eqnarray}\label{eq:stable_111}
\lim_{\tau\downarrow0}\limsup_{N\to\infty} \P[M_{N,\tau}'>\eps
/2]&=&0, \\\label{eq:stable_111a}
\lim_{\tau\downarrow0}\limsup_{N\to\infty}\P[M_{N,\tau}''>\eps/2]&=&0
.
\end{eqnarray}

Let us prove \eqref{eq:stable_111}. Note that the process $\{\mathrm{e}^{\alpha
\eta(t)}, t\geq0\}$ is a martingale. By Doob's maximal
$L^p$-inequality, $\E\|\mathrm{e}^{2\eta}\|_{\infty}\leq C \E \mathrm{e}^{2\eta
(T)}<\infty$.
Thus, $\E\|\mathrm{e}^{\eta_{i,N}}\|_{\infty}^2$ is finite and
\[
\limsup_{N\to\infty} \E M_{N,\tau}'^2 \leq C\lim_{N\to\infty} N
\E\bigl[W_{N}^{2}(0)1_{W_{N}(0)\leq\tau}\bigr]=\frac{C\alpha}{2-\alpha}\tau
^{2-\alpha},
\]
where the last step follows from part \hyperlink{p:1_lem6}{2} of
Proposition \ref{lem:aux_int} with $\kappa=2$.
The right-hand side goes to $0$ as $\tau\downarrow0$. By Chebyshev's
inequality, this proves \eqref{eq:stable_111}.

Let us prove \eqref{eq:stable_111a}. By Theorem \ref{theo:zero}, the
random variable $N^{-1/2}\|\sum_{i=1}^N (\mathrm{e}^{\eta_{i,N}}-\E \mathrm{e}^{\eta
_{i,N}})\|_{\infty}$ converges as $N\to\infty$ to some limiting
(a.s. finite) random variable. Thus, we need to prove that
%
\begin{equation}\label{eq:stable_200}
\lim_{\tau\downarrow0}\limsup_{N\to\infty} \sqrt{N}\E
\bigl[W_{N}(0)1_{W_{N}(0)\leq\tau}\bigr]=0.
\end{equation}
We have, by part \hyperlink{p:1_lem6}{2} of Proposition \ref{lem:aux_int} with
$\kappa=2$,
\[
\limsup_{N\to\infty}N \E\bigl[W_{N}(0)1_{W_{N}(0)\leq\tau}\bigr]^2
\leq
\lim_{N\to\infty}N\E\bigl[W_{N}^2(0)1_{W_{N}(0)\leq\tau}\bigr]=\frac
{\alpha}{2-\alpha}\tau^{2-\alpha}.
\]
This proves \eqref{eq:stable_200} and completes the proof of the lemma.
\end{pf}

\subsection{Completing the proof of the one-sided convergence}
In this section we complete the proof of the one-sided version of
Theorem \ref{theo:stable}.
We will need to introduce some notation. Let $d$ be the Skorokhod
metric on $D[0,T]$. Given a~process~$X$ with sample paths in $D[0,T]$,
we denote by $\mathcal L(X)$ the law of $X$ considered as a~probabili\-ty
measure on $D[0,T]$. Let further $\pi$ be the L\'evy--Prokhorov
distance on the space of~proba\-bility measures on $D[0,T]$. That is,
given two probability measures $\mu_1$ and $\mu_2$ on $D[0,T]$, we define
\[
\pi(\mu_1,\mu_2)=\inf\{\eps>0\dvtx  \mu_1(B)\leq\mu_2(B^{\eps
})+\eps\mbox{ for all Borel } B\subset D[0,T]\},
\]
where $B^{\eps}=\{b\in D[0,T]\dvtx d(b,B)\leq\eps\}$ is the $\eps
$-neighborhood of the set $B$. The next lemma is standard.
\begin{lemma}\label{lem:skor_wspom}
Let $\{X(t), t\in[0,T]\}$ and $\{Y(t), t\in[0,T]\}$ be two
(generally, dependent) stochastic processes with sample paths in
$D[0,T]$, and suppose that for some $\eps>0$,
\[
\P\Bigl[\sup_{t\in[0,T]} |Y(t)|>\eps\Bigr]\leq\eps.
\]
Then, $\pi(\mathcal L(X),\mathcal L(X+Y))\leq\eps$.
\end{lemma}
\begin{pf}
By the definition of the Skorokhod metric, $d(X,X+Y)\leq\sup_{t\in
[0,T]}|Y(t)|$. By assumption, it follows that $\P[d(X,X+Y)>\eps]\leq
\eps$. For every Borel set $B\subset D[0,T]$, we have
\[
\P[X+Y\in B]\leq\P[X\in B^{\eps}]+\P[d(X,X+Y)>\eps]\leq\P[X\in
B^{\eps}]+\eps,
\]
whence the statement of the lemma.
\end{pf}

We are now in position to complete the proof of the one-sided version
of Theorem \ref{theo:stable}, as restated in \eqref
{eq:pr_stab_YN_YYY}. Let $\eps>0$ be fixed. Our aim is to show that
for sufficiently large $N$, we have
%
\begin{equation}\label{eq:pr_stab_wsp0}
\pi(\mathcal L (Y_{N}), \mathcal L(\YYY_{\alpha;\xi}))\leq3\eps.
\end{equation}
By Lemma \ref{lem:est_lower_stat1}, we can find a $\delta>0$ such
that $\P[\sup_{t\in[0,T]}|\YYY_{\alpha;\xi}^{(0,\tau)}(t)|>\eps
]\leq\eps$ for all $\tau<\delta$. By Lemma \ref{lem:skor_wspom}
and \eqref{eq:def_Y_zero_tau}, this implies that for all $\tau<\delta$,
%
\begin{equation}\label{eq:pr_stab_wsp1}
\pi\bigl(\mathcal L\bigl(\YYY_{\alpha;\xi}^{(\tau,\infty)}\bigr), \mathcal
L(\YYY_{\alpha;\xi})\bigr)\leq\eps.
\end{equation}
By Lemma \ref{lem:est_lower_stat2}, we can find $\tau<\delta$ and
$N_1\in\N$ such that $\P[\sup_{t\in[0,T]}|Y_{N}^{(0,\tau
)}(t)|>\eps]\leq\eps$ for $N>N_1$. By Lemma \ref{lem:skor_wspom},
this implies that for all $N>N_1$,
%
\begin{equation}\label{eq:pr_stab_wsp2}
\pi\bigl(\mathcal L\bigl(Y_{N}^{(\tau,\infty)}\bigr), \mathcal{L}(Y_{N})\bigr)\leq\eps.
\end{equation}
By Lemma \ref{lem:est_upper_stat}, we can find $N_1\in\N$ such that
for all $N>N_1$,
%
\begin{equation}\label{eq:pr_stab_wsp3}
\pi\bigl(\mathcal L\bigl(Y^{(\tau,\infty)}_N\bigr), \mathcal L\bigl(\YYY^{(\tau,\infty
)}_{\alpha;\xi}\bigr)\bigr)\leq\eps.
\end{equation}
To complete the proof of \eqref{eq:pr_stab_wsp0}, combine
equations \eqref{eq:pr_stab_wsp1}--\eqref
{eq:pr_stab_wsp3}.

\section*{Acknowledgements}
The author is grateful to Leonid Bogachev for pointing out
reference \cite{cranstonmolchanov05} after \cite{kabluchkoprod} was
completed, and to Ilya Molchanov, Michael Schmutz and Martin Schlather
for useful discussions.

\printhistory

\end{document}